\documentclass[twocolumn]{autart}    % Enable this line and disable the 
\usepackage[utf8]{inputenc}   % codificação UTF-8
\usepackage[T1]{fontenc}      % suporte a caracteres acentuados

\usepackage{amsmath,amssymb,amsfonts,mathtools}
\usepackage{calrsfs,mathrsfs}
\usepackage{graphicx}
\usepackage{xcolor}
\usepackage{float}
\usepackage{multirow}
\usepackage{psfrag}
\usepackage{blindtext}

\usepackage{subcaption}
\newtheorem{assumption}{Assumption}

\newtheorem{theorem}{Theorem}
 
\begin{document}

\begin{frontmatter}

\title{Multivariable Gradient-Based Extremum Seeking Control\\ with Saturation Constraints\thanksref{footnoteinfo}} % Title, preferably not more 

\thanks[footnoteinfo]{This paper was not presented at any IFAC 
meeting.\\ Corresponding author: P.~H.~S.~Coutinho.}

\author[UERJ]{Enzo Ferreira Tomaz Silva}\ead{enzotomazsilva@gmail.com},
\author[UERJ]{Pedro Henrique Silva Coutinho}\ead{phcoutinho@eng.uerj.br}, 
\author[UERJ]{Tiago Roux Oliveira}\ead{tiagoroux@uerj.br},
\author[UCSD]{Miroslav Krsti\'{c}}\ead{mkrstic@ucsd.edu},
\author[LAAS]{Sophie Tarbouriech}\ead{tarbour@laas.fr}.

\address[UERJ]{Dept. of Electronics and Telecommunication Engineering, Rio de Janeiro State University (UERJ), Rio de Janeiro -- RJ, Brazil}
\address[UCSD]{Dept. of Mechanical and Aerospace Engineering, University of California - San Diego (UCSD), La Jolla -- CA, USA}
\address[LAAS]{LAAS-CNRS, Universit\'{e} de Toulouse, CNRS, Toulouse, France}

\begin{keyword}                           % Five to ten keywords,  
Extremum seeking control; Actuator saturation; Gradient algorithm; Multivariable systems; Convex optimization.               % chosen from the IFAC 
\end{keyword}                             % keyword list or with the 
                                          % help of the Automatica 
                                          % keyword wizard

\begin{abstract}                          % Abstract of not more than 200 words.
This paper addresses the multivariable gradient-based extremum seeking control (ESC) subject to saturation. Two distinct saturation scenarios are investigated here: saturation acting on the input of the function to be optimized, which is addressed using an anti-windup compensation strategy, and saturation affecting the gradient estimate. In both cases, the unknown Hessian matrix is represented using a polytopic uncertainty description, and sufficient conditions in the form of linear matrix inequalities (LMIs) are derived to design a stabilizing control gain. The proposed conditions guarantee exponential stability of the origin for the average closed-loop system under saturation constraints. With the proposed design conditions, non-diagonal control gain matrices can be obtained, generalizing conventional ESC designs that typically rely on diagonal structures. Stability and convergence are rigorously proven using the Averaging Theory for dynamical systems with Lipschitz continuous right-hand sides. Numerical simulations illustrate the effectiveness of the proposed ESC algorithms, confirming the convergence even in the presence of saturation.
\end{abstract}

\end{frontmatter}

\section{Introduction}

Extremum seeking control (ESC) is an adaptive, real-time, and model-free optimization strategy. The purpose of this technique is to find an optimal point such that a given desired function (with unknown parameters) is maximized or minimized, that is, its extremum is reached~\cite{LivroTiagoES,ESsurvey}.
Since the first stability analysis for ESC systems provided in~\cite{krstic2000stability}, several efforts have been made to extend the ESC to different classes of maps and control problems, such as time-delay systems~\cite{ArtDelay2}, maps in cascade with partial differential equations~\cite{ArtPDE}, 
% cooperative games with Nash equilibrium~\cite{rodrigues2024sliding,rodrigues2024nash},
and event-triggered control~\cite{rodrigues2025event}. 
% However, these studies deal with unconstrained ESC. 
% do not deal with constraints on actuators.

It is known that, in practice, input constraints can arise due to physical limitations and operational constraints~\cite{tarbouriech2011book,coutinho2020saturation,peixoto2022static}.
If the presence of input constraints is not properly handled in the analysis or synthesis problem, the performance of the closed-loop system might be deteriorated or even lead to unstable behavior.
Within the context of ESC, input constraints have been handled by employing a constrained optimization perspective. In~\cite{frihauf2013finite}, a finite-horizon LQ control problem was solved via ESC by employing the projection operator to handle input constraints and introducing a Newton-based discrete-time ESC. In~\cite{tan2013extremum}, the ESC problem was addressed with a hard saturation nonlinearity constraining the input in a gradient-based ESC scheme for optimizing scalar quadratic maps. Although an anti-windup (AW) compensation is suggested in that work,  the authors claim that it is hard to rigorously demonstrate that the AW mechanism works. 
{To circumvent this issue, the authors proposed penalty-function-based ESC schemes.} 

An AW compensation for ESC has also been employed in an observer-based ESC for a direct-contact membrane distillation process~\cite{eleiwi2018observer}. There, the AW compensation is also exploited in the multivariable case by adding an AW compensator to each input channel, constituting a decentralized AW compensator. For a class of discrete-time nonlinear control systems subject to input constraints, reference~\cite{guay2019extremum} proposes a proportional-integral ESC with a discrete-time AW mechanism. In~\cite{lu2020hardware}, an ESC scheme is proposed for operational control of
mineral grinding, considering both operational indices regulation and throughput maximization. To address the input constraints, a saturation function is applied to each input, and an additional penalty function is introduced to penalize inputs that fall outside the feasible region. A sampled-data ESC for constrained optimization is provided by~\cite{hazeleger2022sampled}. Unlike the aforementioned works that deal with hard saturation constraints, reference \cite{hazeleger2022sampled} deals with the input constraint by employing a barrier function-based method, such that the input constraints are satisfied provided that parametric initialization yields operating conditions that
do not violate the constraints. More recently,~\cite{Karimi2025} investigated the AW penalty-based approach from~\cite{tan2013extremum} for higher derivatives Newton-based ESC schemes under input saturation.

Despite advances in dealing with saturated ESC using AW and penalty-based approaches, existing solutions in the literature are often limited to the scalar case. In particular, when the multivariable case is concerned, a decentralized AW compensation strategy is employed at each input channel.
In addition, none of those papers proposes conditions to design the feedback-adaptation gain and the AW compensation gain. In contrast, some papers, not in the context of ESC, tackle actuator saturation using AW techniques and provide control design conditions in the form of linear matrix inequalities (LMIs)~\cite{GALEANI2009418,Zaccarian2023,Tarbouriech2009,zac:tee/book}. However, these solutions require a deep knowledge of the plant to be analyzed and lack dealing with system optimization, which are strong characteristics of ESC.

%$(i)$ they do not deal with the multivariable case, and $(ii)$ they do not propose control design conditions that ensure the stability of the closed-loop system even in the presence of saturation.

This paper addresses the multivariable gradient-based ESC subject to saturation. Two distinct saturation scenarios are investigated here: saturation acting on the input of the function to be optimized, which is addressed using an AW compensation strategy; and saturation affecting the gradient estimate. %First, for the case where the Hessian matrix is fully unknown, a 
The complete stability analysis is derived for both cases and the convergence of the trajectories to a neighborhood of the optimal point is guaranteed by invoking the averaging theorem for non-differentiable Lipschitz systems~\cite{plotnikov1979averaging}---see also Appendix~\ref{sec:appendix}.  
In particular, the Hessian matrix is considered polytopic uncertain, and control design conditions are established via LMIs to obtain the control gains such that the average system is exponentially stable.
Interestingly, the design methodology presented here allows for non-diagonal control gains, offering to explore greater design flexibility rather than the diagonal gains typically assumed \textit{a priori} in the ESC literature.

\textbf{Notation.}
$\mathbb{R}^n$ is the $n$-dimensional Euclidean space and $\mathbb{R}^{m\times n}$ the set of real matrices of order $m \times n$. A symmetric positive (negative) definite matrix $X$ is denoted by $X > 0~(X < 0)$. For a matrix $X \in \mathbb{R}^{n \times m}$, $X^\top$ denotes its transpose, and $X_{(\ell)} \in \mathbb{R}^{1 \times m}$ denotes its $\ell$-th row, and $X_{(i,j)}$ the element in the $i$-th row and $j$-th column. For a vector $x \in \mathbb{R}^{n}$, $x_\ell \in \mathbb{R}$ denotes the $\ell$-th element of $x$.

\section{Extremum Seeking Control under Actuator Saturation}
\label{sec:preliminaries}

\begin{figure}[!ht]
\begin{center}
\includegraphics[width=\columnwidth]{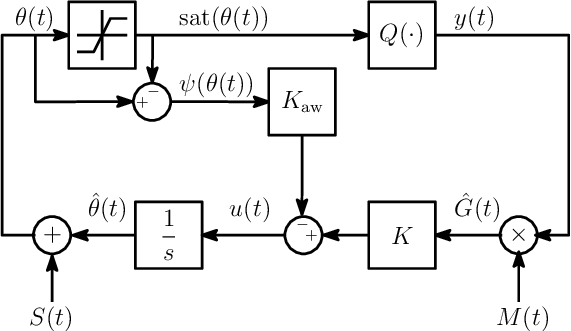}
\caption{Extremum seeking control system with saturation in the input map.} 
\label{fig:diag_ESC_input_constrained}
\end{center}
\end{figure}

Consider the multivariable gradient-based ESC system with input saturation shown in Fig.~\ref{fig:diag_ESC_input_constrained}. 

We deal with the following nonlinear static map under input saturation constraints:
\begin{align} \label{eq:map_static_input_constrained}
       &y(t) = Q(\mathrm{sat}(\theta(t))  \\
       &= Q^* + \frac{1}{2}(\mathrm{sat}(\theta(t)) - \theta^*)^\top H(\mathrm{sat}(\theta(t)) - \theta^*), \label{VitoriaVascao}
\end{align}
where $Q^* \in \mathbb{R}$ is the unknown optimal point of the map, 
$\theta^* \in \mathbb{R}^n$ is the unknown optimizer of the map, 
$\theta \in \mathbb{R}^n$ is the input vector, 
$H \in \mathbb{R}^{n \times n}$ is the unknown Hessian matrix, 
$y \in \mathbb{R}$ is the output of the map, and $\mathrm{sat}(\cdot) : \mathbb{R}^n \to \mathbb{R}^n$ is the element-wise saturation function defined as follows: 
\begin{align}\label{eq:sat_function}
    \mathrm{sat}(\theta) = \begin{bmatrix}
        \mathrm{sat}(\theta_1) \\ \vdots \\ \mathrm{sat}(\theta_n)
    \end{bmatrix} = 
    \begin{bmatrix}
        \mathrm{sign}(\theta_1) \min(|\theta_1|,\overline{\theta}_1) \\ \vdots \\ \mathrm{sign}(\theta_n) \min(|\theta_n|,\overline{\theta}_n)
    \end{bmatrix},
\end{align}
where $\overline{\theta}_\ell > 0$ is the saturation bound of the 
$\ell$-th input signal.

\begin{assumption} \label{assump:theta*}
The unknown optimizer $\theta_\ell^*$ lies within the {region in which saturation does not occur, defined by}
\begin{align} \label{eq:theta_barra_&_theta*}
    {\Theta^* = \{\theta^* \in \mathbb{R}^n : |\theta^*_{\ell}| < \overline{\theta}_{\ell}, \, \ell = 1,\ldots,n\}.}
\end{align}
\end{assumption}
%
%{In particular, equation (\ref{VitoriaVascao}) states that %if maps with saturated inputs are at least locally quadratic, our %methodology can be applied and can provide guarantees within some %neighborhood of the extremum. This requirement is not restrictive, since %any nonlinear function that is twice continuously differentiable can be %approximated by a quadratic function near its extremum. Therefore, all %stability results derived in this paper hold at least locally. For maps %that are not locally quadratic and may not yield exponential stability of %the average system, the approach in \cite{TAN2006889} could lead to %asymptotic practical stability (rather than exponential), but the form of %averaging theory used therein is not available for systems on Banach spaces %(e.g., those involving saturation functions).} 
%
{ In particular, equation (\ref{VitoriaVascao}) shows that if maps with saturated inputs are locally quadratic, our methodology applies and provides guarantees in a neighborhood of the extremum. This assumption is mild, as any twice continuously differentiable nonlinear function $Q(\cdot)$ admits a quadratic approximation near its extremum. Hence, all stability results in this paper hold at least locally. For maps that are not locally quadratic and may not yield exponential stability of the average system, the approach in \cite{TAN2006889} could ensure asymptotic practical (rather than exponential) stability, but its averaging framework is not available for Banach-space systems (e.g., those with saturation functions).}
{Notice also that for a twice continuously differentiable nonlinear function $Q(\cdot)$, the Hessian matrix is symmetric.}

\subsection{Probing and Demodulation Signals}

In the scheme shown in Fig~\ref{fig:diag_ESC_input_constrained}, the probing and demodulation dithers are respectively defined by~\cite{ghaffari2012multivariable}:
\begin{align} 
    S(t) = \begin{bmatrix}
       a_1\sin{(\omega_1 t)} & \cdots & a_n\sin{(\omega_n t)} \\
    \end{bmatrix}^\top, \label{eq:vetor_S_input_contrained} \\
    M(t) = \begin{bmatrix}
       \frac{2}{a_1}\sin{(\omega_1 t)} & \cdots & \frac{2}{a_n}\sin{(\omega_n t)} \\
    \end{bmatrix}^\top, \label{eq:vetor_M_input_constrained}
\end{align}
where $a_i$, for all $i = 1,\ldots,n$, are non-zero amplitudes of the dither signals, and their frequencies are selected according to the following assumption.

\begin{assumption}\label{assump:probing_w}
    For a given angular frequency $\omega>0$, the probing frequencies are selected such that 
\begin{align} \label{eq:relacao_sinais_omega_pert_input}
    \omega_i = \omega_i' \omega, \quad i = 1,\ldots,n,
\end{align}
where $\omega_i'$ is a rational number that satisfies
    \begin{align}\label{eq:probing_w}
        \omega_i' \notin \left\lbrace\omega_j', \frac{1}{2}(\omega_j'+\omega_k'), \omega_j' + 2 \omega_k', \omega_k'\pm \omega_l' \right\rbrace,
    \end{align}
for all $i,j,k = 1,\ldots,n$ and $l$.
\end{assumption}
{To guarantee convergence in multivariable ESC schemes, it is necessary to choose a sufficiently large $\omega$ as well as distinct probing frequencies ($\omega_i \neq \omega_j$), and ensure that the ratio $\omega_i/\omega_j$ is rational and $\omega_i + \omega_j \neq \omega_k$ for different $i, j, k$. These conditions are fulfilled under Assumption~\ref{assump:probing_w}, according to~\cite{ghaffari2012multivariable}.}

Under Assumption~\ref{assump:probing_w}, the gradient estimation $\hat{G}(t)$ is driven by the periodic perturbations and demodulation strategy so that
\begin{align}\label{eq:gradient_1}
    \hat{G}(t) = M(t) y(t).
\end{align}

\subsection{\textcolor{black}{Extremum Seeking with Anti-Windup Design}}

The output of the integrator provides the estimate $\hat{\theta}(t) \in \mathbb{R}^n$ of the optimum point $\theta^*$, such that the estimation error is given by
\begin{align}\label{eq:estimation_error}
    \tilde{\theta}(t) = \hat{\theta}(t) - \theta^*.
\end{align}
Then, the estimation error dynamics is
\begin{align} \label{eq:theta_til_ponto_input_constrained}
    \dot{\tilde{\theta}}(t) = \dot{\hat{\theta}}(t) = u(t).
\end{align}
Different from the unconstrained ESC, an AW compensation term is introduced into the input of the integrator, such that the AW compensator is given by 
\begin{align}\label{eq:control_law_1}
    u(t) = K \hat{G}(t) - K_{\mathrm{aw}} \psi(\theta(t)),
\end{align}
where $K \in \mathbb{R}^{n \times n}$ is the feedback gain, which is not assumed diagonal, as usually done in ESC works, and $K_{\mathrm{aw}} \in \mathbb{R}^{n \times n}$ is the AW compensation gain, and 
\begin{equation} \label{eq:def_sat_input}
    \psi(\theta(t)) = \theta(t) - \mathrm{sat}(\theta(t))
\end{equation}
is the dead-zone nonlinearity~\cite{tarbouriech2011book}. 
The introduction of the AW compensation aims to drive the output of the quadratic map $y(t)$ to the optimal point $Q^*$ when the input operates in the constrained region. Notice that $\psi(\theta) \equiv 0$ when the input operates in the {linear region defined in~\eqref{eq:theta_barra_&_theta*}}, and~\eqref{eq:control_law_1} reduces to a standard ESC law.

\subsection{Reshaping the Gradient Estimate}

Using the relation
\begin{equation} \label{eq:state_map}
    \theta(t) = \hat{\theta}(t) + S(t),
\end{equation}
together with~\eqref{eq:vetor_S_input_contrained},~\eqref{eq:def_sat_input}, and~\eqref{eq:estimation_error}, it follows that 
\begin{equation}\label{eq:signals_relation}
        \mathrm{sat}(\theta(t)) - \theta^* =  S(t) + \tilde{\theta}(t) - \psi(\theta(t)). 
\end{equation}

By substituting \eqref{eq:map_static_input_constrained} into~\eqref{eq:gradient_1} and using~\eqref{eq:signals_relation}, we obtain
% \begin{align}\label{eq:gradient_2}
%     &\hat{G}(t) = M(t) Q^* + \frac{1}{2}M(t)S^\top(t)HS(t) \nonumber \\
%     &+ \frac{1}{2}M(t)S^\top(t)H\tilde{\theta}(t) 
%     - \frac{1}{2}M(t)S^\top(t)H\psi(\theta(t)) \nonumber \\
%     &+ \frac{1}{2}M(t)\tilde{\theta}^\top(t)HS(t) 
%     + \frac{1}{2}M(t)\tilde{\theta}^\top(t)H\tilde{\theta}(t) \nonumber \\
%     &- \frac{1}{2}M(t)\tilde{\theta}^\top(t)H\psi(\theta(t)) 
%     - \frac{1}{2}M(t)\psi^\top(\theta(t)) HS(t) \nonumber \\
%     &- \frac{1}{2}M(t)\psi^\top(\theta(t)) H\tilde{\theta}(t) + \frac{1}{2}M(t)\psi^\top(\theta(t)) H\psi(\theta(t)) 
% \end{align}
\begin{align}\label{eq:gradient_3}
    &\hat{G}(t) = M(t) Q^* + \frac{1}{2}M(t)S^\top(t)HS(t) \nonumber \\
    &+ M(t)S^\top(t)H\tilde{\theta}(t) 
    - M(t)S^\top(t)H\psi(\theta(t)) \nonumber \\
    &+ \frac{1}{2}M(t)\tilde{\theta}^\top(t)H\tilde{\theta}(t) \nonumber - M(t)\tilde{\theta}^\top(t)H\psi(\theta(t)) \nonumber \\
    & + \frac{1}{2}M(t)\psi^\top(\theta(t)) H\psi(\theta(t)). 
\end{align}
Let $\Omega(t) \coloneqq M(t)S^\top(t)H$. By noticing that~\cite{rodrigues2025event}
\begin{equation} \label{eq:delta_input_constrained}
   \Omega(t) = M(t)S^\top(t)H = H + \Delta(t)H,
\end{equation}
where $\Delta(t)$ is a matrix-valued function whose entries are
\begin{align}
    \Delta_{ii}(t) &= 1-\cos(2\omega_it), \\
    \Delta_{ij}(t) &= {\frac{a_{j}}{a_{i}}\cos((\omega_i-\omega_j)t) - \frac{a_{j}}{a_{i}}\cos((\omega_i+\omega_j)t)}, \label{coesmeedamiao}
\end{align}
we can rewrite~\eqref{eq:gradient_3} as follows:
% \begin{align}\label{eq:gradient_4}
%     &\hat{G}(t) = M(t) Q^* + \frac{1}{2}M(t)S^\top(t)HS(t) \nonumber \\
%     &+ (H+\Delta(t)H)\tilde{\theta}(t) 
%     - (H+\Delta(t)H)\psi(\theta(t)) \nonumber \\
%     &+ \frac{1}{2}M(t)\tilde{\theta}^\top(t)H\tilde{\theta}(t) \nonumber - M(t)\tilde{\theta}^\top(t)H\psi(\theta(t)) \nonumber \\
%     & + \frac{1}{2}M(t)\psi^\top(\theta(t)) H\psi(\theta(t)). 
% \end{align}
\begin{align}\label{eq:gradient_5}
    &\hat{G}(t) = (H+\Delta(t)H)\tilde{\theta}(t) 
    - (H+\Delta(t)H)\psi(\theta(t)) + w(t), 
\end{align}
where
\begin{align}
    &w(t) {=} M(t) Q^* {+} \frac{1}{2}{\Omega(t)} S(t) {+} \frac{1}{2}M(t)\tilde{\theta}^\top(t)H\tilde{\theta}(t) \nonumber \\
    &- M(t)\tilde{\theta}^\top(t)H\psi(\theta(t)) 
    + \frac{1}{2}M(t)\psi^\top(\theta(t)) H\psi(\theta(t)).
\end{align}
% Since the term~$\tilde{\theta}^\top(t)H\tilde{\theta}(t)$ is quadratic in $\tilde{\theta}(t)$, it can be neglected in a local analysis~\cite{ariyur2003real}. 
Then, the gradient estimate~\eqref{eq:gradient_5} can be still rewritten as
\begin{align}\label{eq:gradient_6}
    \hat{G}(t) &= H \tilde{\theta}(t)  - H \psi(\theta(t)) \nonumber \\
    & + \Delta(t)H \tilde{\theta}(t) - \Delta(t)H \psi(\theta(t)) + w(t). 
\end{align}

\subsection{Closed-Loop System}

By substituting~\eqref{eq:control_law_1} and~\eqref{eq:gradient_6} into~\eqref{eq:theta_til_ponto_input_constrained}, we obtain the following closed-loop dynamics:
\begin{align}\label{eq:closed_loop_1}
    \dot{\tilde{\theta}}(t) &= K H \tilde{\theta}(t) - K H \psi(\theta(t))  - K_{\mathrm{aw}} \psi(\theta(t)) \nonumber \\
    & + K \Delta(t)H \tilde{\theta}(t) - K \Delta(t)H \psi(\theta(t)) + K w(t). 
\end{align}
%The stability analysis of~\eqref{eq:gradient_5} can now be analyzed using the averaging method. To this purpose, it is necessary to express (\ref{eq:closed_loop_1}) in an appropriate form by means of a new time-scale to evaluate the effect of $\omega$ in the dynamics. It is also important to note that the time-varying disturbances $w(t)$ and $\Delta(t)$ both have zero mean values, which will enable the subsequent averaging analysis.
{The stability analysis of~(\ref{eq:closed_loop_1}) can now be performed using the averaging method. To this purpose, it is necessary to express in an appropriate form by means of a new time-scale to evaluate the effect of $\omega$ in the dynamics. It is also important to note that the time-varying disturbances $w(t)$ and $\Delta(t)$ both have zero mean values, which will enable the subsequent averaging analysis.}

\subsection{Rescaling of Time} \label{subsec:Defining_new_time}

For the stability analysis of the closed-loop system, a change in the time scale is performed. According to~\eqref{eq:relacao_sinais_omega_pert_input}, it is ensured that the dither frequencies and their combinations are rational. 
Thus, there exists a period
\begin{align}
    T = 2\pi \times \mathrm{LCM}\left\lbrace\frac{1}{\omega_1}, \ldots, \frac{1}{\omega_n}\right\rbrace,
\end{align}
where $\mathrm{LCM}$ denotes the least common multiple. The change of time scale of the system in \eqref{eq:closed_loop_1} consists of a transformation %$\tau = \omega t$, 
\begin{equation} \label{tau_duro}
\tau = \omega t\,,
\end{equation}
where
\begin{align}\label{eq:angular_frequency}
    \omega := \frac{2\pi}{T}.
\end{align}
Hence, reminding that $\theta(t)=\tilde{\theta}(t)+ S(t) + \theta^*$ from (\ref{eq:estimation_error}) and (\ref{eq:state_map}), the right-hand side of \eqref{eq:closed_loop_1} can be rewritten as a function of $\tilde{\theta}$:  %
\begin{align} \label{eq:esc_temp_mudanca_Theta_til}
    \frac{d\tilde{\theta}\left(\tau\right)}{d\tau} = \frac{1}{\omega}\mathcal{F}\left(\tau, \tilde{\theta}, \frac{1}{\omega}\right), 
\end{align}
where 
\begin{align}\label{eq:esc_temp_mudanca_Theta_til_function}
      &\mathcal{F}\left(\tau,  \tilde{\theta}, \frac{1}{\omega}\right)   = K H \tilde{\theta}(\tau) - K H \psi(\theta(\tau))  - K_{\mathrm{aw}} \psi(\theta(\tau)) \nonumber \\
    & + K \Delta(\tau)H \tilde{\theta}(\tau) - K \Delta(\tau)H \psi(\theta(\tau)) + K w(\tau).
\end{align}

\subsection{Average Closed-Loop System}

After performing the time scaling, the average version of \eqref{eq:esc_temp_mudanca_Theta_til}--\eqref{eq:esc_temp_mudanca_Theta_til_function} can be computed as follows 
\begin{align} \label{eq:Sistema_Medio_Derivado_input}
    \frac{d\tilde{\theta}_\mathrm{av}(\tau)}{d\tau} &= \frac{1}{\omega}\mathcal{F}_\mathrm{av}(\tilde{\theta}_\mathrm{av}), \\
    \mathcal{F}_\mathrm{av}(\tilde{\theta}_\mathrm{av}) &= \frac{1}{T} \int_0^T   \mathcal{F}(\delta, \tilde{\theta}_\mathrm{av},0) d\delta.  \label{eq:Media_calculada_Sistema_input}
\end{align}
% For each term, the average is computed below:
% \begin{align} 
%     S_\mathrm{av}(\tau) &= \frac{1}{T} \int_0^T   S(\delta) d\delta = 0, \\  \dot{S}_\mathrm{av}(\tau) &= \frac{1}{T} \int_0^T   \dot{S}(\delta) d\delta = 0, \label{eq:Media_Sinal_S_input} \\
%    M_\mathrm{av}(\tau) &= \frac{1}{T} \int_0^T   M(\delta) d\delta = 0, \\
%    \dot{M}_\mathrm{av}(\tau) &= \frac{1}{T} \int_0^T   \dot{M}(\delta) d\delta = 0, \label{eq:Media_Sinal_M_input} \\
%     \Delta_\mathrm{av}(\tau) &= \frac{1}{T} \int_0^T  \Delta(\delta) d\delta = 0, \\
%     \dot{\Delta}_\mathrm{av}(\tau) &= \frac{1}{T} \int_0^T   \dot{\Delta}(\delta) d\delta = 0. \label{eq:Media_Sinal_Delta_input}
% \end{align}
% As a result, one can obtain
% \begin{align*}
%     \Omega_{\mathrm{av}}(\tau) &= \frac{1}{T} \int_0^T \Omega(\delta) d\delta = H,  \\
%     \dot{\Omega}_{\mathrm{av}}(\tau) &= \frac{1}{T} \int_0^T \dot{\Omega}(\delta) d\delta = 0.
% \end{align*}
{For a sufficiently large $\omega > 0$, we can ``freeze'' the average state of $\tilde{\theta}(t)$ and treat it as constant in~\eqref{eq:Media_calculada_Sistema_input}. Then, from Assumption~\ref{assump:probing_w}, we have that
\begin{align}
    \frac{1}{T} \int_0^T \Delta_{ij}(t) dt = 0, \quad
    \frac{1}{T} \int_0^T w_i(t)dt  = 0,
\end{align}
for all $i,j = 1,\ldots,n$. Therefore, we obtain the following average dynamics for the closed-loop system~\eqref{eq:closed_loop_1}: 
% \begin{small}
\begin{equation} \label{eq:theta_til_medio}
     \!\!\dot{\tilde{\theta}}_\mathrm{av}(t) \!=\! KH\tilde{\theta}_\mathrm{av}(t) {-} KH\psi(\theta_\mathrm{av}(t)) {-} K_\mathrm{aw}\psi(\theta_\mathrm{av}(t)),
\end{equation}
where $\theta_\mathrm{av}(t) = \tilde{\theta}_\mathrm{av}(t) + \theta^*$, provided that $S(t)$ has zero mean over a period $T$.}
% \end{small}

{The aim of deriving the average version of the closed-loop system is to design the control gains $K$ and $K_{\mathrm{aw}}$ such that the linear time-invariant system under actuator saturation~\eqref{eq:theta_til_medio} is globally exponentially stable, and then investigate the local stability analysis of the non-autonomous time-varying system~\eqref{eq:closed_loop_1} using the Averaging Theory for systems with Lipschitz continuous right-hand sides~(see Appendix~\ref{sec:appendix}).}

{An interesting aspect to be noticed in~\eqref{eq:theta_til_medio} is that simply taking some diagonal matrix $K$ --- as done in the context free of constraints \cite{ghaffari2012multivariable} --- such that the matrix $K H$ is Hurwitz does not necessarily ensure the exponential stability of the origin of the average closed-loop system~\eqref{eq:theta_til_medio}. Thus, this aspect emphasizes the importance of developing a constructive method to design the control gains $K$ and~$K_{\mathrm{aw}}$.}

\subsection{Polytopic Embedding of the Hessian Matrix}

In general, ESC methodologies rely on the assumption that the unknown Hessian matrix $H$ is a negative (or positive) definite matrix, depending on whether the optimal point is a maximum (or a minimum). Based on this, a negative (or positive definite) diagonal structure is assigned to the gain matrix $K$. However, this approach is not straightforward for the input-constrained case, especially because the AW gain $K_{\mathrm{aw}}$ should be designed as well.

As an alternative to designing the control gains $K$ and $K_{\mathrm{aw}}$, we propose here to exploit a polytopic embedding of the Hessian matrix, as stated in the following Assumption. 
\begin{assumption}\label{assump:H*}
    The unknown Hessian matrix $H$ takes values in the following polytopic domain:
    \begin{align}\label{eq:polytpoic_H_input}
        \mathcal{H} = \mathrm{co}\{H_1, \ldots,H_N\},
    \end{align}
    where $N$ is the number of vertices of the polytope and $H_i$, $i = 1,\ldots,N$, are known matrices.
\end{assumption}
Under Assumption~\ref{assump:H*}, any unknown Hessian matrix $H \in \mathcal{H}$ can be parameterized as follows:
\begin{align}\label{eq:H_parameterization}
    H =  H(\alpha) = \sum_{i=1}^N \alpha_i H_i,
\end{align}
where $\alpha = (\alpha_1,\ldots,\alpha_N)$ is the vector of fixed but unknown parameters that belongs to the unitary simplex
\begin{align}\label{eq:unity_simplex}
    \Xi = \left\lbrace \alpha \in \mathbb{R}^N : \sum_{i = 1}^N \alpha_i = 1, \; \alpha_i \geq 0, i = 1,\ldots,N \right\rbrace.
\end{align}
% and $H_i \in \mathbb{R}^{n \times n}$, for all $i = 1,\ldots, N$, are the polytope vertices, which are known matrices. Thus, it is possible to obtain the following uncertain polytopic description for the closed-loop system:
% \begin{multline} \label{eq:average_system_politopo_input}
%      \dot{\tilde{\theta}}_\mathrm{av}(\tau) = \frac{1}{\omega}KH(\alpha)\tilde{\theta}_\mathrm{av}(\tau) - \frac{1}{\omega}KH(\alpha)\psi(\theta_\mathrm{av}(\tau))\\ -\frac{1}{\omega}K_\mathrm{aw}\psi(\theta_\mathrm{av}(\tau)).
% \end{multline}
With the polytopic parameterization in~\eqref{eq:H_parameterization}, the Hessian matrix is still unknown; only the vertices of the polytopic domain need to be known. The Hessian polytope can be obtained using different uncertainty representations. 
For instance, if one assumes that $\lambda_1 I \leq H \leq \lambda_2 I$,
we can simply select $H_1 = \lambda_1 I$ and $H_2 = \lambda_2 I$, which results in a polytope with two vertices. Another strategy is to assign an uncertainty to an estimate of the nominal Hessian $H_0$, such that $H = H_0 + \delta H_0$, $|\delta| \leq \overline{\delta}$, which can be embedded in a two-vertex polytope with $H_1 = (1-\overline{\delta})H_0$ and $H_2 = (1+\overline{\delta})H_0$. 
We can also assign different bounds to the elements of an uncertain Hessian matrix. For instance, if $H$ is of order two as 
\begin{align}\label{eq:H_example}
    H = \begin{bmatrix}
    h_{11} & h_{12} \\ h_{12} & h_{22}
\end{bmatrix}
\end{align}
and we assume that $|h_{11}| \leq \overline{\delta}_1$, $|h_{12}| \leq \overline{\delta}_2$, and $|h_{22}| \leq \overline{\delta}_3$, we can construct a polytope with eight vertices ($2^p$, where $p=3$ is the number of uncertain parameters, $h_{11}$, $h_{12}$, and $h_{22}$) from the combination of the bounds of those uncertain parameters.
In a more general sense, if we consider an affine representation as
\begin{align}\label{eq:H_affine}
    H = \Gamma_{0} + \delta_1 \Gamma_{1} + \ldots + \delta_p \Gamma_{p},
\end{align}
where $|\delta_i| \leq \overline{\delta}_i$ and $\Gamma_i$, $i = 0,1,\ldots,p$, are known matrices, we can obtain a polytopic representation for $H$ with $2^p$ vertices.
For instance, the uncertain matrix~\eqref{eq:H_example} can be rewritten in the affine form~\eqref{eq:H_affine} with $\delta_1 = h_{11}$, $\delta_2 = h_{12}$, $\delta_3 = h_{22}$, 
\begin{align*}
    \Gamma_0 = 
    \begin{bmatrix}
        0 & 0 \\ 0 & 0
    \end{bmatrix}, \;
    \Gamma_1 = 
    \begin{bmatrix}
        1 & 0 \\ 0 & 0
    \end{bmatrix}, \;
    \Gamma_2 = 
    \begin{bmatrix}
        0 & 1 \\ 1 & 0
    \end{bmatrix}, \;
    \Gamma_3 = 
    \begin{bmatrix}
        0 & 0 \\ 0 & 1
    \end{bmatrix}.     
\end{align*}
Another approach for estimating the Hessian matrix is using an averaging- or perturbation-based approach via dither signals from~\cite{ghaffari2012multivariable}. However, this method depends on the unconstrained feedback operation to estimate the Hessian matrix. Thus, this method might not be directly employed in the presence of saturation.

\subsection{Stability Analysis}

In this section, we provide stabilization conditions to 
design the control gains $K \in \mathbb{R}^{n \times n}$ and $K_{\mathrm{aw}} \in \mathbb{R}^{n \times n}$, such that the origin of the average closed-loop system~\eqref{eq:theta_til_medio} is {globally} exponentially stable. Then, by invoking the Averaging Theory (see Appendix~\ref{sec:appendix}), we show the {asymptotic convergence of the trajectories of the ESC system under input saturation~\eqref{eq:closed_loop_1} to a neighborhood of the extremum point.}

%%%%%%%%%%%%%%%%%%%%%%%%%%%%%%%%%%%%%%%%%%%%%%%%%%%

\subsubsection{Stabilization of the Average Closed-Loop System}
\label{sec:Main_Results_Input_Constrained}

% In this section, a stability criterion for the averaged system is derived using a sector condition for the dead-zone nonlinearity $\psi(\cdot)$. Thereafter, by applying the averaging theorem from~\cite{plotnikov1979averaging}, a convergence condition for the original system is established. Finally, the control design task is expressed as an LMI-based optimization problem.
%This subsection introduces a stability condition for the average system based on a sector condition for the dead-zone function. After that, it is established a convergence condition of the original system using the averaging theorem of~\cite{plotnikov1979averaging}. Finally, the control design condition is formulated as an optimization problem based on LMIs.

% \subsubsection{Stability analysis of the average system}

To establish the condition to design the control gains $K$ and $K_{\mathrm{aw}}$ in~\eqref{eq:theta_til_medio}, we exploit the generalized sector condition~\cite{tarbouriech2011book} of the dead-zone nonlinearity $\psi(\theta(t))$.
% \begin{lem}[Generalized sector condition~\cite{tarbouriech2011book}]
%     If $v$ and $w$ are elements of 
%     \begin{align}
%         \mathcal{S}(v-w,\overline{u}) = 
%         \{v, w \in \mathbb{R}^m : |v-w| \leq \overline{u}\},
%     \end{align}
%     then the nonlinearity $\psi(v)$ satisfies the following inequality:
%     \begin{align}
%         \psi^\top(v) T \left(\psi(v) - w\right)
%     \end{align}
% \end{lem}
Based on this result, the following lemma {proposes a global sector condition}. %is established.
\begin{lem}\label{lem:input_constrained}
    If Assumption~\ref{assump:theta*} holds, then
    % then $\theta_{\mathrm{av}}$ and $\tilde{\theta}_{\mathrm{av}}$ are elements of the following set
    % \begin{multline} \label{eq:conjunto_THETA}
    %             \Theta = \lbrace \theta_{\mathrm{av}}, \tilde{\theta}_{\mathrm{av}} \in \mathbb{R}^{n} : \ |\theta_{\mathrm{av}(\ell)}-\tilde{\theta}_{\mathrm{av}(\ell)}| \leq \overline{\theta}_{\ell},\\ \, \ell = 1,\ldots,n\rbrace,
    % \end{multline}
    % which ensures that the sector condition 
    \begin{align}\label{eq:setor_psi_input_constrained}
        \psi^\top(\theta_{\mathrm{av}}) \Lambda \left(\psi(\theta_{\mathrm{av}}) -  \tilde{\theta}_{\mathrm{av}}\right) \leq 0,
    \end{align}
    holds for any diagonal positive definite matrix~$\Lambda \in \mathbb{R}^{n \times n}$, and any $\theta_{\mathrm{av}}, \tilde{\theta}_{\mathrm{av}} \in \mathbb{R}^n$.
\end{lem}
\begin{pf}
    {
    Under Assumption~\ref{assump:theta*}, we have $-\overline{\theta}_\ell \leq \theta_\ell^* \leq \overline{\theta}_\ell$. By evaluating~\eqref{eq:signals_relation} in the average sense, we have 
    $\theta_{\mathrm{av}} - \tilde{\theta}_{\mathrm{av}} = \theta^*$, provided that $S(t)$ has zero mean over a period~$T$ as in \eqref{eq:angular_frequency}. Thus, Assumption~\ref{assump:theta*} ensures that  
    $\theta_{\mathrm{av}}$ and $\tilde{\theta}_{\mathrm{av}}$ are elements of the following set
    \begin{multline} \label{eq:conjunto_THETA}
                \Theta = \lbrace \theta_{\mathrm{av}}, \tilde{\theta}_{\mathrm{av}} \in \mathbb{R}^{n} : \ |\theta_{\mathrm{av}(\ell)}-\tilde{\theta}_{\mathrm{av}(\ell)}| \leq \overline{\theta}_{\ell},\\ \, \ell = 1,\ldots,n\rbrace,
    \end{multline}    
    which implies that
    \begin{align} \label{eq:saturation_interval}
        -\overline{\theta}_{\ell} \leq \theta_{\mathrm{av}(\ell)}-\tilde{\theta}_{\mathrm{av}(\ell)} \leq \overline{\theta}_{\ell}.
    \end{align}}
    %The rest of the proof follows similar steps to~\cite[Lemma~1.6]{tarbouriech2011book} and is omitted here for brevity. ~\hfill$\blacksquare$
    %Consider now the following three cases.
    The following three cases are now taken into account.
    %Para mostrar que~\eqref{eq:setor_psi} é satisfeita, três casos serão analisados:
    \begin{itemize}
        \item \textbf{Case 1}: $\theta_{{\mathrm{av}}{(\ell)}} > \overline{\theta}_{\ell}$. The following holds
    \begin{align}
        \psi(\theta_{{\mathrm{av}}{(\ell)}}) = \theta_{{\mathrm{av}}{(\ell)}} - \overline{\theta}_{\ell} > 0.
    \end{align}        
     It follows from~\eqref{eq:saturation_interval} that
    \begin{align}
         \psi(\theta_{{\mathrm{av}}{(\ell)}})\! - \!\tilde{\theta}_{{\mathrm{av}}{(\ell)}} \!=\! \theta_{{\mathrm{av}}{(\ell)}}-\!\tilde{\theta}_{{\mathrm{av}}{(\ell)}}\!  -\! \overline{\theta}_{\ell}\! \leq \!0.
    \end{align}
    Thus, $ \psi^\top(\theta_{{\mathrm{av}(\ell)}}) \Lambda_{(\ell,\ell)} \left(\psi(\theta_{{\mathrm{av}(\ell)}}) -  \tilde{\theta}_{{\mathrm{av}(\ell)}}\right) \leq 0$, provided that $\Lambda_{(\ell,\ell)} > 0$.\\
        \item \textbf{Case 2}: $ -\overline{\theta}_{\ell} \leq \theta_{{\mathrm{av}}{(\ell)}} \leq \overline{\theta}_{\ell}$, the dead-zone function $\psi(\theta_{{\mathrm{av}}{(\ell)}})$ is zero, since $\mathrm{sat}(\theta_{{\mathrm{av}}{(\ell)}}) = \theta_{{\mathrm{av}}{(\ell)}}$. In this case, we obtain 
        \begin{align*}
          \psi^\top\!(\!\theta_{{\mathrm{av}(\ell)}}\!) \Lambda_{(\ell,\ell)}\! \left(\!\psi(\theta_{{\mathrm{av}(\ell)}})\! -\!  \tilde{\theta}_{{\mathrm{av}(\ell)}}\!\right)\! = \!0, \; \forall\Lambda_{(\ell,\ell)}.  
        \end{align*}
        \item \textbf{Case 3}: $\theta_{{\mathrm{av}}{(\ell)}} <- \overline{\theta}_{\ell}$. The following holds
    \begin{align}
        \psi(\theta_{{\mathrm{av}}{(\ell)}}) = \theta_{{\mathrm{av}}{(\ell)}} + \overline{\theta}_{\ell} < 0.
    \end{align}        
     It follows from~\eqref{eq:saturation_interval} that
    \begin{align}
         \psi(\theta_{{\mathrm{av}}{(\ell)}})\! - \!\tilde{\theta}_{{\mathrm{av}}{(\ell)}} \!=\! \theta_{{\mathrm{av}}{(\ell)}}-\!\tilde{\theta}_{{\mathrm{av}}{(\ell)}}\!  +\! \overline{\theta}_{\ell}\! \geq \!0.
    \end{align}
    Thus, $ \psi^\top(\theta_{{\mathrm{av}(\ell)}}) \Lambda_{(\ell,\ell)} \left(\psi(\theta_{{\mathrm{av}(\ell)}}) -  \tilde{\theta}_{{\mathrm{av}(\ell)}}\right) \leq 0$, provided that $\Lambda_{(\ell,\ell)} > 0$.\\
    \end{itemize}
    From the three cases presented, it can be verified that the inequality in~\eqref{eq:setor_psi_input_constrained} is satisfied for all $\theta_{\mathrm{av}}$ and $\tilde{\theta}_{\mathrm{av}}$ in~\eqref{eq:conjunto_THETA}, which is provided since Assumption~\ref{assump:theta*} holds. This concludes the proof.~\hfill$\blacksquare$
\end{pf}

An interesting aspect of the result stated in Lemma~\ref{lem:input_constrained} is the sector condition of $\psi(\theta_{\mathrm{av}})$ established in~\eqref{eq:setor_psi_input_constrained} with respect to $\tilde{\theta}_{\mathrm{av}}$, which is the variable of the dynamics under study in~\eqref{eq:theta_til_medio}.

Based on the result established in Lemma~\ref{lem:input_constrained}, we provide in the sequel a stabilization condition to design the gains of the control law~\eqref{eq:control_law_1} to ensure the exponential stability of the average closed-loop system~\eqref{eq:theta_til_medio} in the presence of saturation with the AW compensation.

\begin{lem}\label{thm:1_input_constrained}
    Consider the closed-loop input constrained ESC system in~\eqref{eq:closed_loop_1} under Assumptions~\ref{assump:theta*},~\ref{assump:probing_w},~and~\ref{assump:H*}. 
    Given a positive scalar $\eta > 0$, if there exist a symmetric positive definite matrix  $P \in \mathbb{R}^{n \times n}$, a diagonal positive definite matrix $\Lambda \in \mathbb{R}^{n \times n}$, and matrices $Z, Z_{\mathrm{aw}} \in \mathbb{R}^{n \times n}$, such that
    % a control gain $K \in \mathbb{R}^{n \times n}$ and matrices $P = P^\top > 0 \in \mathbb{R}^{n \times n}$,
    % $Q = Q^\top > 0 \in \mathbb{R}^{n \times n}$, there are matrices
    % $L \in \mathbb{R}^{n \times n}$, a diagonal matrix $T > 0 \in \mathbb{R}^{n \times n}$ and $K_\mathrm{aw} \in \mathbb{R}^{n \times n}$ such that
    % \begin{align}\label{eq:ineq-stability_input}
    %     \begin{bmatrix}
    %         {H}^\top K^\top P + P K H + 2 \eta P & \star\\
    %         T L - PK_\mathrm{aw} - PKH  & - 2 T
    %     \end{bmatrix} < 0,
    % \end{align}    
    % \begin{equation}
    %  \begin{bmatrix}\label{eq:input_constrained_LMI}
    %    P     &  0 &  L_{(\ell)}^\top \\
    %    \star & Q & I_{(\ell)^\top} \\
    %    \star & \star &\bar \theta_{(\ell)}^{2}                                 
    %  \end{bmatrix} \geq  0, \quad \ell  = 1,2,\ldots,n,
    % \end{equation}
    \begin{equation}\label{eq:LMI_thm_input}
    \begin{bmatrix} 
         ZH_i+ H_i Z^\top + 2\eta P  & \star \\
         \Lambda - Z_{\mathrm{aw}}^\top - {H_i Z^\top}  & -2\Lambda
    \end{bmatrix}<0, \, \forall i = 1,\ldots,N,
    \end{equation}
    % \begin{align}
    %     \begin{bmatrix}
    %             P & 0 & Y_{(\ell)}^\top \\
    %             \star & Q & T_{(\ell)}^\top \\
    %             \star & \star & \overline{\theta}_\ell^2(2T - I)
    %     \end{bmatrix} \geq 0, \quad \ell = 1,\ldots,n.  
    %     \label{eq:saturation_thm_input} 
    % \end{align}
    then, the origin of the average closed-loop system~\eqref{eq:theta_til_medio}, with $K = P^{-1}Z$ and $K_\mathrm{aw} = P^{-1}Z_{\mathrm{aw}}$, is {globally} exponentially stable with decay rate $\eta$, that is:
    \begin{align}\label{eq:theta_til_medio_inequecao}
        \|\tilde{\theta}_{\mathrm{av}}(t)\| \leq \kappa e^{-\eta t}\|\tilde{\theta}_{\mathrm{av}}(0)\|,
    \end{align}
    where $\kappa = \sqrt{\lambda_{\max}(P) / \lambda_{\min}(P)}$. 
    % provided that $\tilde{\theta}_{\mathrm{av}}(0)$ is taken inside the region $\mathcal{E}_\theta \subset \Theta$, with
    % \begin{align}
    %     \mathcal{E}_\theta = \{\tilde{\theta}_{\mathrm{av}} \in \mathbb{R}^n : V(\tilde{\theta}_{\mathrm{av}}) \leq 1 - \theta_{\mathrm{av}}^\top(0) Q \theta_{\mathrm{av}}(0)\},
    % \end{align}
    % where 
    % \begin{align} \label{eq:lyapunov_Theta_av}
    %     V(\tilde{\theta}_{\mathrm{av}}) = \tilde{\theta}_{\mathrm{av}}^\top P \tilde{\theta}_{\mathrm{av}}
    % \end{align}
    % is a Lyapunov function that certifies the exponential stability of the origin of the average closed-loop system~\eqref{eq:theta_til_medio}.
\end{lem}
\begin{pf}
Assume that conditions~\eqref{eq:LMI_thm_input} hold. Provided that $\alpha \in \Xi$, with $\Xi$ given in~\eqref{eq:unity_simplex}, it follows from~\eqref{eq:LMI_thm_input} and Assumption~\ref{assump:H*} that
\begin{align}\label{eq:thm1-proof-1}
   \begin{bmatrix}
        ZH+ {H} Z^\top + 2\eta P  & \star \\
        \Lambda - Z_{\mathrm{aw}}^\top - H Z^\top   & -2 \Lambda
   \end{bmatrix} < 0.
\end{align}
By substituting $Z = P K$, $Z_{\mathrm{aw}} = P K_{\mathrm{aw}}$ in~\eqref{eq:thm1-proof-1}, we obtain
\begin{align}\label{eq:thm1-proof-2}
   \begin{bmatrix}
         P K H+ {H} K^\top P + 2\eta P  & \star \\
        \Lambda - K_{\mathrm{aw}}^\top P - H K^\top P   & -2\Lambda
   \end{bmatrix} < 0.
\end{align}
By multiplying~\eqref{eq:thm1-proof-2} with $[\tilde{\theta}_{\mathrm{av}}^\top \; \psi^\top(\theta_{\mathrm{av}})]$ on the left and its transpose on the right, it results in
\begin{align}
     &\left[KH \tilde{\theta}_\mathrm{av} - (K_\mathrm{aw} + KH)\psi(\theta_\mathrm{av}))\right]^{\top} P \tilde{\theta}_{\mathrm{av}} \nonumber \\ 
     &+ \tilde{\theta}_{\mathrm{av}}^{\top} P \left[KH\tilde{\theta}_\mathrm{av} - (K_\mathrm{aw} + KH)\psi(\theta_\mathrm{av})\right] \nonumber  \\ 
     &+ 2 \eta \tilde{\theta}_{\mathrm{av}}^{\top} P \tilde{\theta}_{\mathrm{av}} - 2\psi^\top(\theta_{\mathrm{av}}) \Lambda \left(\psi(\theta_{\mathrm{av}}) - \tilde{\theta}_{\mathrm{av}}\right) < 0.
\end{align}
Thus, provided that $\tilde{\theta}_{\mathrm{av}}$ and ${\theta}_{\mathrm{av}}$ are elements of $\Theta$ in~\eqref{eq:conjunto_THETA}, under Assumption~\ref{assump:theta*}, it follows from Lemma~\ref{lem:input_constrained} that  
\begin{align}\label{eq:exp-ineq_input}
    \dot{V}(\tilde{\theta}_{\mathrm{av}}(t)) + 2 \eta V(\tilde{\theta}_{\mathrm{av}}(t)) < 0,
\end{align}
where 
% {$V(\tilde{\theta}_{\mathrm{av}})$} is given in~\eqref{eq:lyapunov_Theta_av}.
{
\begin{align} \label{eq:lyapunov_Theta_av}
    V(\tilde{\theta}_{\mathrm{av}}) = \tilde{\theta}_{\mathrm{av}}^\top P \tilde{\theta}_{\mathrm{av}}
\end{align}
is a Lyapunov function that certifies the exponential stability of the origin of the average closed-loop system~\eqref{eq:theta_til_medio}.}
From the Comparison Lemma \cite{KH:02}, it follows from~\eqref{eq:exp-ineq_input} that
\begin{align}
    {V}(\tilde{\theta}_{\mathrm{av}}(t)) \leq e^{-2 \eta t} {V}(\tilde{\theta}_{\mathrm{av}}(0)). 
\end{align}
Furthermore, since 
\begin{align}
    \lambda_{\min}(P) \|\tilde{\theta}_{\mathrm{av}}\|^2 \leq 
    V(\tilde{\theta}_{\mathrm{av}}) \leq \lambda_{\max}(P) \|\tilde{\theta}_{\mathrm{av}}\|^2,
\end{align}
it is possible to show that~\eqref{eq:theta_til_medio_inequecao} holds.
% \begin{align} \label{eq:theta_til_medio_inequecao}
%     \|\tilde{\theta}_{\mathrm{av}}(t)\| \leq \kappa e^{-\eta t }\|\tilde{\theta}_{\mathrm{av}}(0)\|,
% \end{align}
% where $\kappa = \sqrt{\lambda_{\max}(P) / \lambda_{\min}(P)}$.
Then, the origin of the system is exponentially stable. This concludes the proof.~\hfill$\blacksquare$
% Moreover, using the fact that $(I-T)^\top (I-T) \geq 0$, we have
% $T^2 \geq 2 T - I$. Using this property, it follows from~\eqref{eq:saturation_thm_input} that
% \begin{align}
%     \begin{bmatrix}
%         P & 0 & Y_{(\ell)}^\top \\
%         \star & Q & T_{(\ell)}^\top \\
%         \star & \star & \overline{\theta}_\ell^2 T^2
%     \end{bmatrix} \geq 0, \quad \ell = 1,\ldots,n.  
%     \label{eq:thm1-proof-3} 
% \end{align}
% Then, using the change of variables $Y = LT$ and multiplying~\eqref{eq:thm1-proof-3} with $\mathrm{diag}(I,$
\end{pf}

\subsubsection{Asymptotic Convergence to a Neighborhood of the Extremum}

Lemma~\ref{thm:1_input_constrained} established a condition to design the control gains that render the origin of the average closed-loop system~\eqref{eq:theta_til_medio} exponentially stable. In the sequel, we state the main result of this section, which provides the local asymptotic convergence to a neighborhood of the extremum by employing the Averaging Theory (see Appendix~\ref{sec:appendix}).

\begin{theorem}\label{thm:main_result_1}
    Consider the ESC system in Fig.~\ref{fig:diag_ESC_input_constrained} with locally quadratic nonlinear map (\ref{eq:map_static_input_constrained})--(\ref{VitoriaVascao}) subject to input saturation and  
    %Consider
    the corresponding average closed-loop dynamics~\eqref{eq:theta_til_medio} %of the ESC system %subject to input %saturation~\eqref{eq:theta_til_medio}
    under Assumptions~\ref{assump:theta*},~\ref{assump:probing_w}, and~\ref{assump:H*}. 
    If the conditions of Lemma~\ref{thm:1_input_constrained} 
    %condition (\ref{eq:ineq-stability_input}) and Assumptions~\ref{assump:theta*} to \ref{assump:H*}
    are all satisfied,  
    then, for $\omega > 0$ sufficiently large in (\ref{eq:relacao_sinais_omega_pert_input}) and $a_i>0$ sufficiently small in (\ref{eq:vetor_S_input_contrained})--(\ref{eq:vetor_M_input_constrained}), there exist constants $\eta > 0$ and $\kappa$ in (\ref{eq:theta_til_medio_inequecao}), such that: 
    % the equilibrium $\tilde{\theta}_{\mathrm{av}} = 0$ is {globally} exponentially stable. %and $\tilde{\theta}_{\mathrm{av}}(t)$ converges exponentially to zero.
    % In particular, there exist constants \textcolor{black}{$\kappa, \eta > 0$} such that: 
    \begin{align} \label{eq:theta_desigualdades_input}
        \|\theta(t) - \theta^\ast\| \leq {\kappa} e^{-\eta t} \|{\theta}(0)-\theta^*\| + \mathcal{O}\left(a + \frac{1}{\omega}\right), \\
        \label{eq:y_desigualdades_input}
        \lim_{t \to \infty}\sup\textcolor{black}{|y(t) - Q^*| =  \mathcal{O}\left(a^2 + \frac{1}{\omega^2}\right)},
        %|y(t) - Q^\ast| \leq  \kappa_y e^{-2\eta t} + \mathcal{O}\left(a^2 + \frac{1}{\omega^2}\right),
    \end{align}
    with $a = \sqrt{\sum_{i=1}^n a_i^2}$.  
%    
    %$\kappa_\theta$ and \textcolor{red}{$\kappa_y$ %constants} which depends of the initial condition %$\theta(0)$. 
\end{theorem}
\begin{pf}  
    %Rewriting the Lyapunov function in ~\eqref{eq:lyapunov_Theta_av} as
    %\begin{align}
    %    V(\tilde{\theta}_{\mathrm{av}}) = \tilde{\theta}_{\mathrm{av}}^\top P \tilde{\theta}_{\mathrm{av}},
    %\end{align}
   % where $P$ is a symmetric and positive definite matrix. Thus, it is possible to find
   % \begin{align} \label{ref_thetal_til_converge_zero_input}
   %     \|\tilde{\theta}_{\mathrm{av}}(\tau)\| \leq \kappa_\theta e^{-\eta \tau / \omega} \|\tilde{\theta}_{\mathrm{av}}(0)\|,
   % \end{align}
   % 
    Since the differential equation in~\eqref{eq:esc_temp_mudanca_Theta_til} has Lipschitz continuous right-hand sides, %discontinuity on the right side,
    due to the presence of the saturating function, and the closed-loop average system \eqref{eq:theta_til_medio} is exponentially stable from Lemma~\ref{thm:1_input_constrained}, %and Lipschitz continuous,
    by applying averaging theorem in \cite{plotnikov1979averaging} (see also Appendix~\ref{sec:appendix}, with $\varepsilon:=1/\omega$), it follows that: 
    \begin{align} \label{diff_theta_theta_avg}
        \textcolor{black}{\|\tilde{\theta}(t) - \tilde{\theta}_{\mathrm{av}}(t)\| \leq \mathcal{O}\left(\frac {1}{\omega}\right) , \quad \forall t\geq 0,}
    \end{align}
    for $\omega$ sufficiently large. 
%%
%    \begin{align} \label{diff_theta_theta_avg}
%        \textcolor{red}{\|{\theta}(t) - {\theta}_{\mathrm{av}}(t)\| \leq %\mathcal{O}\left(a+ \frac {1}{\omega}\right).}
%    \end{align}
    Then, applying the triangle inequality into (\ref{diff_theta_theta_avg}), {from the relation~\eqref{eq:theta_til_medio_inequecao}}, we can obtain 
    \begin{align} \label{eq:demonstracao_estabilidade_theta_input}
        \| \tilde{\theta}(t) \| \leq \kappa  e^{-\eta t} \|\tilde{\theta}(0)\| + \mathcal{O}\left(\frac{1}{\omega}\right).
    \end{align}
From~\eqref{eq:state_map} and the definition of $\tilde{\theta}(t)$ in (\ref{eq:estimation_error}), we can write the following equivalence 
    \begin{align} \label{relacao_theta_til&theta_medio}
        \theta(t) - \theta^* = \tilde{\theta}(t)+S(t),
    \end{align}
to obtain (\ref{eq:theta_desigualdades_input}) since $S(t)$ in (\ref{eq:vetor_S_input_contrained}) is of order $\mathcal{O}(a)$, with $a = \sqrt{\sum_{i=1}^n a_i^2}$.

Now, consider (\ref{eq:map_static_input_constrained}) to write the following output error 
\begin{align}
    \tilde{y}(t) := y(t) - Q^*,~~~~y(t) = Q( \mathrm{sat}(\theta(t))).
\end{align}
By computing its norm, and using the Cauchy–Schwarz inequality, one gets
\begin{multline}
   |\tilde{y}(t)| = |( \mathrm{sat}(\theta(t)) - \theta^*)^\top H(\mathrm{sat}(\theta(t)) - \theta^*)| \\ \leq \|H\| \| \mathrm{sat}(\theta(t)) - \theta^*\|^2.
\end{multline} 
Using the dead-zone function definition in \eqref{eq:def_sat_input}, we obtain 
\begin{align} \label{eq:output_error}
     |\tilde{y}(t)| \leq \|H\| \| \theta(t) -  \psi(\theta(t))  - \theta^*\|^2.
\end{align}
From~\eqref{eq:theta_desigualdades_input},  inequality~\eqref{eq:output_error} can be reformulated as
\begin{multline} \label{eq:lim_output_error}
     \lim_{t \to \infty} \sup|\tilde{y}(t)| \leq \lim_{t \to \infty} \sup \|H\| \| \theta(t)  - \theta^* -  \psi(\theta(t)) \|^2 \leq \\ 
     \lim_{t \to \infty} \sup \|H\| \left[ \left\| \mathcal{O}\left(\!a \!+\! \frac{1}{\omega}\right)\right\|^2 \!\!+\!  \left\|\psi\!\left(\!\mathcal{O}\left(a \!+\! \frac{1}{\omega}\right) \!+\! \theta^*\right) \right\|^2\right].
\end{multline}
Using Assumption~\ref{assump:theta*}, the dead-zone function $$\psi\left(\mathcal{O}\left(a + \frac{1}{\omega}\right) + \theta^*\right)=0,$$ for $\omega$ sufficiently large and $a$ sufficiently small, due to the condition~\eqref{eq:theta_barra_&_theta*}, with {$|\theta^*_{\ell}| \leq |\mathcal{O}\left(a + \frac{1}{\omega}\right)| + |\theta_{\ell}^*| < \bar{\theta}_{\ell}$}. Hence, by employing Young's inequality to the term $ \left\| \mathcal{O}\left(a \!+\! \frac{1}{\omega}\right)\right\|^2$, we obtain 
\begin{equation}
    \lim_{t \to \infty}  \sup |\tilde{y}(t)| = \mathcal{O}\left(a^2 + \frac{1}{\omega^2}\right),
\end{equation}
leading to (\ref{eq:y_desigualdades_input}), which   
%As a result, the inequalities \eqref{eq:theta_desigualdades_input} and %\eqref{eq:y_desigualdades_input} are guaranteed.
completes the proof.~\hfill$\blacksquare$
\end{pf}

%%%%% Artigo com a Saturação no Gradiente %%%%%
\section{\textcolor{black}{Extremum Seeking Control under Gradient Saturation}}
\label{sec:3}

\begin{figure}[!ht]
\begin{center}
\includegraphics[width=8cm]{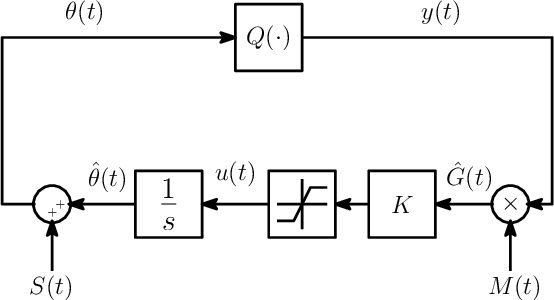}
\caption{Extremum seeking control system \textcolor{black}{under gradient saturation.}} 
\label{fig:diag_ES_sat}
\end{center}
\end{figure}

Consider the multivariable gradient-based ESC under \textcolor{black}{gradient saturation} shown in Fig.~\ref{fig:diag_ES_sat}.

{In this feedback system, we are assuming the unknown multi-input nonlinear map is locally quadratic in the neighborhood of its extremum, such that} 
\begin{equation} \label{eq:map_static}
       y(t) = Q(\theta(t)) =  Q^* + \frac{1}{2}(\theta(t) - \theta^*)^\top H(\theta(t) - \theta^*),
\end{equation}
where $Q^*$, $\theta^*$, and $H$ of the map are defined as in~\eqref{eq:map_static_input_constrained}. 
%where $Q^* \in \mathbb{R}^n$ is the unknown optimal point of the map, 
%$\theta^* \in \mathbb{R}^n$ is the unknown optimizer of the map, 
%$\theta \in \mathbb{R}^n$ is the input vector, 
%$H \in \mathbb{R}^{n \times n}$ is the unknown Hessian Matrix and 
%$y \in \mathbb{R}^n$ is the map output.
%Even though the Hessian matrix is unknown, it can be assumed that it is a positive definite matrix when the minimum point is desired and a negative definite matrix when the maximum point is desired.
%
Also, the input vector $\theta(t)$ applied to the multivariable static map is defined as in~\eqref{eq:state_map} and the estimation error $\tilde{\theta}(t)$ is given by \eqref{eq:estimation_error}. We also consider the probing and demodulation signals $S(t)$ and $M(t)$ according to \eqref{eq:vetor_S_input_contrained} and \eqref{eq:vetor_M_input_constrained},  respectively, satisfying Assumption~\ref{assump:probing_w}. 
%\begin{equation} \label{eq:input_map}
%    \theta(t) = \hat{\theta}(t) + S(t)
%\end{equation}
%where $\hat{\theta} \in \mathbb{R}^n$ is the estimated value of $\theta^*$.

\subsection{\textcolor{black}{Extremum Seeking with Bounded Update Rates}}

Different from Section~\ref{sec:preliminaries}, where we studied the saturation in the input of the map, in this section, we deal with saturation in the \textcolor{black}{gradient estimate}
\begin{align}\label{eq:control_law_1_cosmeedamiao}
    %u(t) = K \hat{G}(t).
    u(t) = \mathrm{sat}(K\hat{G}(t)).
\end{align}
In this case, the dynamics for $\tilde{\theta}(t)$ is described as follows:
\begin{align} \label{eq:theta_til_ponto_&_saturacao}
   %\dot{\tilde{\theta}}(t) = \dot{\hat{\theta}}(t) = \mathrm{sat}(u(t)) = \mathrm{sat}(K\hat{G}(t)),
   \dot{\tilde{\theta}}(t) = \dot{\hat{\theta}}(t) = u(t) = \mathrm{sat}(K\hat{G}(t)),
\end{align}
where $K \in \mathbb{R}^{n \times n}$ is the control gain to be designed, $\hat{G}(t)$ is given as in~\eqref{eq:gradient_1}, and $\mathrm{sat}(\cdot)$ is the saturation function defined in the element-wise sense as in~\eqref{eq:sat_function}:
 \begin{align}
     \mathrm{sat}(K \hat{G}) \!\!=\!\! \begin{bmatrix}
         \mathrm{sat}(K \hat{G}_1) \\ \vdots \\ \mathrm{sat}(K \hat{G}_n)
     \end{bmatrix} \!\!=\!\! 
     \begin{bmatrix}
         \mathrm{sign}(K \hat{G}_1) \min(|K \hat{G}_1|,\overline{u}_1) \\ \vdots \\ \mathrm{sign}(K \hat{G}_n) \min(|K \hat{G}_n|,\overline{u}_n)
     \end{bmatrix}\!,
 \end{align}
with $\overline{u}_\ell > 0$ being the limit of the $\ell$-th \textcolor{black}{gradient estimate.} 
In contrast to Section~\ref{sec:preliminaries}, this section considers the presence of saturation in the \textcolor{black}{gradient estimate}. Note that dealing with saturation before the integration is also a important problem, as it leads to a feedback loop with bounded update rates. Indeed, the primary motivation for considering gradient saturation lies in establishing the first formulation of classical ESC with bounded update rates. In this setting, saturation should not be regarded as a limitation but rather as a beneficial mechanism that enables such boundedness, such as in the Lie-Bracket bounded ESC approaches \cite{DURR20131538}, \cite{SCHEINKER201425}, \cite[Ch.~6]{SCHEINKER201425_book}. 

According to the gradient estimation expression $\hat{G}(t)=M(t)y(t)$ with the quadratic map~\eqref{eq:map_static} and using the relation $\theta(t) - \theta^* = \tilde{\theta}(t)+S(t)$, we have that
% As $y(t)$ is the quadratic map output, the gradient estimate can be written as:
\begin{equation}
    \hat{G}(t) = M(t)\left(Q^* + \frac{1}{2}(\tilde{\theta}(t) + S(t))^\top H(\tilde{\theta}(t) + S(t))\right),
\end{equation}
or still
\begin{multline} \label{eq:grad_simp}
     \hat{G}(t) = M(t)Q^* + \frac{1}{2}M(t)\tilde{\theta}^\top(t)H\tilde{\theta}(t)\\
     + M(t)S^\top(t)H\tilde{\theta}(t) + \frac{1}{2}M(t)S^\top(t)HS(t). 
\end{multline}
Analogously, using the matrices $\Omega(t)$ and $\Delta(t)$ defined according to~\eqref{eq:delta_input_constrained}--\eqref{coesmeedamiao},  
%By defining the matrix
%\begin{equation} \label{eq:delta}
%    \Omega(t) = M(t)S^\top(t)H,
%\end{equation}
%the multiplication in \eqref{eq:delta} results in a matrix of the following %form:
%\begin{equation} \label{eq:delta_final}
%    \Omega(t) = H + \Delta(t)H,
%\end{equation}
%where
%    $\Delta_{ii} = 1-\cos(2\omega_it)$,
%    $\Delta_{ij} = \frac{a_{j}}{a_{i}}\cos(\omega_i-\omega_j) - \frac{a_{j}}{a_{i}}\cos(\omega_i+\omega_j)$.
%
equation \eqref{eq:grad_simp} is expressed as 
\begin{multline} \label{eq:grad_simp_2}
     \hat{G}(t) = M(t)Q^* + \frac{1}{2}M(t)\tilde{\theta}^\top(t) H\tilde{\theta}(t) \\
     +  \Omega(t)\tilde{\theta}(t) + \frac{1}{2}\Omega(t)S(t).
\end{multline}
Since the term~$\tilde{\theta}^\top(t)H\tilde{\theta}(t)$ is quadratic in $\tilde{\theta}(t)$, it can be neglected in a local analysis~\cite{ariyur2003real}.
Then, the dynamics of \eqref{eq:grad_simp_2} can be rewritten as:
\begin{equation} \label{eq:grad_derivado}
     \dot{\hat{G}}(t) = H\mathrm{sat}(K\hat{G}(t)) + \Delta(t)H\mathrm{sat}(K\hat{G}(t)) + \varsigma(t),
\end{equation}
where
\begin{multline} \label{eq:resto_grad_derivado}
        \varsigma(t) = \dot{M}(t)Q^* + \dot{\Delta}(t)H\tilde{\theta}(t) + \frac{1}{2}H\dot{S}(t)+ \frac{1}{2} \dot{\Delta}(t)HS(t)  \\ + \frac{1}{2}\Delta(t)H\dot{S}(t).
\end{multline}

\subsection{Defining a New Time Scale for Averaging}
By adopting an analogous procedure as in  Section~\ref{subsec:Defining_new_time} and noticing that $\varsigma(t)$ in~\eqref{eq:resto_grad_derivado} has zero mean over a period $T:=2\pi/\omega$ given as in~\eqref{eq:angular_frequency}, the following average dynamics is obtained for the new time scale $\tau=\omega t$ from~\eqref{eq:grad_derivado}:
\begin{equation} \label{eq:grad_medio}
     \dot{\hat{G}}_{\mathrm{av}}(\tau) = \frac{1}{\omega} H u_{\mathrm{av}}(\tau) = \frac{1}{\omega}H\mathrm{sat}(K\hat{G}_{\mathrm{av}}(\tau)).
\end{equation}
Based on the dead-zone parametrization for the saturation nonlinearity discussed in~\cite{tarbouriech2011book}, (\ref{eq:control_law_1_cosmeedamiao}) can be written in terms of 
\begin{equation} \label{eq:def_sat}
    \psi(K\hat{G}) = K\hat{G} - \mathrm{sat}(K\hat{G}),
\end{equation}
and the average closed-loop system obtained in~\eqref{eq:grad_medio} can be rewritten as
\begin{align}\label{eq:average_system1}
    \dot{\hat{G}}_{\mathrm{av}}(\tau) = \frac{1}{\omega}H K \hat{G}_{\mathrm{av}}(\tau) - \frac{1}{\omega}H \psi(K \hat{G}_{\mathrm{av}}(\tau)),
\end{align}
where $u_{\mathrm{av}} = \mathrm{sat}(K \hat{G}_{\mathrm{av}})$.

Similarly to Section~\ref{sec:preliminaries}, the Hessian matrix is also assumed to satisfy Assumption~\ref{assump:H*} and 
%In general, solutions available in the literature are developed for the stability analysis of extremum seeking control systems, assuming the knowledge of the sign of the Hessian matrix $H$. Based on this, a diagonal structure with the opposite sign is assigned to the gain matrix $K$. Although this approach requires little knowledge of the Hessian matrix $H$, it becomes difficult to design the gain matrix using constructive design techniques via LMIs. 
%For this purpose, it is assumed that the Hessian matrix $H$ is unknown, but takes values within a polytopic set according to the following parameterization:
%\begin{align}\label{eq:polytpoic_H}
%    H = H(\alpha) = \sum_{i=1}^N \alpha_i H_i,
%\end{align}
%where the vector of uncertain parameters $\alpha = (\alpha_1,\ldots,\alpha_N)$ belongs to the unitary simplex
%\begin{align}
%    \Lambda = \left\lbrace \alpha \in \mathbb{R}^N : \sum_{i = 1}^N \alpha_i = 1, \; \alpha_i \geq 0, i = 1,\ldots,N \right\rbrace
%\end{align}
%and $H_i \in \mathbb{R}^{n \times n}$, $i = 1,\ldots, N$ are the polytope vertices, that are known matrices. 
%
% Thus, it is possible to obtain
the following uncertain polytopic description for the average closed-loop system is finally obtained:
 \begin{align} \label{eq:dinamica_Gav_politopo}
     % \dot{\hat{G}}_{\mathrm{av}}(\tau) = \frac{1}{\omega}H(\alpha) K \hat{G}_{\mathrm{av}}(\tau) - \frac{1}{\omega}H(\alpha) \psi(K \hat{G}_{\mathrm{av}}(\tau))
     {\dot{\hat{G}}_{\mathrm{av}}(t) = H(\alpha) K \hat{G}_{\mathrm{av}}(t) - H(\alpha) \psi(K \hat{G}_{\mathrm{av}}(t)),}
 \end{align}
{where $H(\alpha)$ satisfies the Assumption~\ref{assump:H*} and the parameterization given in~\eqref{eq:H_parameterization}--\eqref{eq:unity_simplex}.} 
% The objective here is to determine a robust control gain $K \in \mathbb{R}^{n \times n}$ guaranteeing exponential stability of the average closed-loop dynamics. Subsequently, by applying the averaging theorem, the stability of the corresponding non-average closed-loop system~\eqref{eq:grad_derivado} is proved.

%The problem addressed in this paper is to design a control robust gain $K \in \mathbb{R}^{n \times n}$ such that the closed-loop average system is exponentially stable. Then, by invoking the averaging theorem, we prove the stability of the non-averaged closed-loop system~\eqref{eq:grad_derivado}.

%%%%%%%%%%%%%%%%%%%%%%%%%%%%%%%%%%%%%%%%%%%%%%%%%%%

\subsection{Stability Analysis}
\label{sec:Main_Results}

In this section, we provide a stabilization condition to 
design the control gain $K \in \mathbb{R}^{n \times n}$ such that the origin of the average closed-loop system~\eqref{eq:average_system1}, or equivalently (\ref{eq:dinamica_Gav_politopo}), is exponentially stable. Then, by invoking the Averaging Theory (see Appendix~\ref{sec:appendix}), we show that the trajectories of the ESC system \textcolor{black}{under gradient saturation} converge exponentially to a neighborhood of the extremum point.

% This section presents a stability condition for the average system based on a sector condition for the dead-zone function. After that, a convergence condition of the original system is established using the averaging theorem of~\cite{plotnikov1979averaging}. Finally, the control design condition is formulated as an optimization problem based on LMIs.

\subsubsection{Stabilization of the Average Closed-Loop System}

%The following lemma establishes a sector condition for the dead-zone nonlinearity
The following lemma {revisits the sector condition in~\cite[Lemma~1.6]{tarbouriech2011book} for} the dead-zone nonlinearity $\psi(K \hat{G})$.
\begin{lem}\label{lem:saturacao}
    Consider a matrix $L \in \mathbb{R}^{m \times n}$. If $\hat{G}_{\mathrm{av}}$ is an element of the set
    \begin{align}\label{eq:conjunto_G}
        \mathcal{G} = \left\lbrace \hat{G}_{\mathrm{av}} \in \mathbb{R}^n : |(K-L)_{(\ell)} \hat{G}_{\mathrm{av}}| \leq \overline{u}_{\ell}, \, \ell = 1,\ldots,n \right\rbrace,
    \end{align}
    then
    \begin{align}\label{eq:setor_psi}
        \psi^\top(K \hat{G}_{\mathrm{av}}) \Upsilon \left(\psi(K \hat{G}_{\mathrm{av}}) - L \hat{G}_{\mathrm{av}}\right) \leq 0,
    \end{align}
    for any diagonal positive definite matrix $\Upsilon \in \mathbb{R}^{n \times n}$.
\end{lem}
\begin{pf}
    {The proof follows similar steps to~\cite[Lemma~1.6]{tarbouriech2011book} and is omitted here for brevity.}~\hfill$\blacksquare$
\end{pf}
%
% \subsubsection{Control Design Conditions}
With Lemma~\ref{lem:saturacao}, we develop the stabilization condition to design the control gain $K$ that renders the origin of the average closed-loop system~\eqref{eq:average_system1} exponentially stable {in a regional context}. This stabilization condition is stated in the following Lemma.

%The theorem below provides a constructive LMI-based condition for designing the gain of the extremum-seeking control system.
\begin{lem}\label{lem:Gav_stability}
Consider the average closed-loop dynamics of the ESC system \textcolor{black}{under gradient saturation}~\eqref{eq:average_system1} under Assumptions~\ref{assump:probing_w} and~\ref{assump:H*}.
{Let $\eta > 0$ and $\epsilon > 0$ be given scalars.}
If there exist a symmetric positive definite matrix $W \in \mathbb{R}^{n \times n}$, a diagonal positive definite matrix $\widetilde{\Upsilon} \in \mathbb{R}^{n \times n}$, and matrices $X, Y, Z \in \mathbb{R}^{n \times n}$, such that the following inequalities hold:
\begin{equation}\label{eq:LMI_thm}
    \begin{bmatrix} 
         H_i Z + Z^\top H_i + 2\eta W & \star & \star\\
         W-X^{\top}+\epsilon H_i Z & -\epsilon (X^{\top}+X) & \star \\
         Y - \widetilde{\Upsilon} H_i & -\epsilon \widetilde{\Upsilon} H_i & -2 \widetilde{\Upsilon}
    \end{bmatrix}<0,
\end{equation}
for all $i = 1,\ldots,N$ and
\begin{align}
    \begin{bmatrix}
            W &  Z_{(\ell)}^\top - Y_{(\ell)}^\top \\
            \star & \overline{u}_\ell^2
    \end{bmatrix} \geq 0, \quad \ell = 1,\ldots,n,  
    \label{eq:saturation_thm} 
\end{align}
then, the origin of the average closed-loop system~\eqref{eq:average_system1} with $K = Z X^{-1}$ is exponentially stable and the region 
\begin{align}\label{eq:elipse_V}
   \mathcal{E} = \{\hat{G}_{\mathrm{av}} \in \mathbb{R}^n : V(\hat{G}_{\mathrm{av}}) \leq 1 \},
\end{align}
is a subset of $\mathcal{G}$ in~\eqref{eq:conjunto_G} {with $L = Y X^{-1}$}, where
\begin{equation}\label{eq:lyapunov_Gav}
    V(\hat{G}_{\mathrm{av}}) = \hat{G}_{\mathrm{av}}^{\top} P \hat{G}_{\mathrm{av}},
\end{equation}
with $P = X^{-\top} W X^{-1}$, is a Lyapunov function that certifies the exponential stability of the origin of~\eqref{eq:average_system1}.
Thus, any trajectory $\hat{G}_{\mathrm{av}}(t)$ with initial condition $\hat{G}_{\mathrm{av}}(0) \in \mathcal{E}$ satisfy 
\begin{align}
    \|\hat{G}_{\mathrm{av}}(t)\| \leq \kappa_g e^{-\eta t}\|\hat{G}_{\mathrm{av}}(0)\|,
\end{align}
where $\kappa_g = \sqrt{\lambda_{\max}(P) / \lambda_{\min}(P)}$.~\hfill$\blacksquare$
%-------------------Antiga Condição----------------------%
% and $\mathcal{E} \subset \mathcal{G} \cap \mathcal{D}$,where
%\begin{align}
%    \mathcal{D} = \{ \hat{G}_{\text{av}} \in \mathbb{R}^n : h_j^\top \hat{G}_{\text{av}} \leq 1, \ldots, n_f \} 
%\end{align}
%is an polyhedral set which contains the origin.
%-------------------Antiga Condição----------------------%
\end{lem}
\begin{pf}
%----------Nova condição para o controlador----------%
Assume that conditions~\eqref{eq:LMI_thm} and~\eqref{eq:saturation_thm} hold.
Provided that $\alpha \in \Xi$, with $\Xi$ given in~\eqref{eq:unity_simplex}, it follows from~\eqref{eq:LMI_thm} and Assumption~\ref{assump:H*} that
\begin{align}\label{eq:proof-lemma4-1}
    \begin{bmatrix} 
         H Z + Z^\top H + 2\eta W & \star & \star\\
         W-X^{\top}+\epsilon H Z & -\epsilon (X^{\top}+X) & \star \\
         Y - \widetilde{\Upsilon} H & -\epsilon \widetilde{\Upsilon} H & -2 \widetilde{\Upsilon}
    \end{bmatrix}<0.
\end{align}
From~\eqref{eq:LMI_thm}, we have that $X + X^\top > 0$, which ensures that $X$
is invertible. Moreover, as $\widetilde{\Upsilon}$ is positive definite, it is also invertible. It allows us to multiply~\eqref{eq:proof-lemma4-1} by $\mathrm{diag}(X^{-\top}, X^{-\top}, \widetilde{\Upsilon}^{-1})$ on the left and its transpose on the right, which results in
% \begin{small}
\begin{align}\label{eq:proof-lemma4-2}
    \begin{bmatrix} 
         \Psi_{11} & \star & \star\\
         \Psi_{21} & -\epsilon (X^{-\top}+X^{-1}) & \star \\
         \Upsilon L - H X^{-1} & -\epsilon H X^{-1}  & -2 \Upsilon
    \end{bmatrix}<0,
\end{align}
% \end{small}
$\!\!$where $\Psi_{11} = X^{-\top} H K {+} K^\top H X^{-1} {+} 2\eta P$, $\Psi_{21} = P-X^{-1} + \epsilon X^{-\top} H K$, $K = Z X^{-1}$, $P = X^{-\top}W X^{-1}$, $L = Y X^{-1}$, and $\Upsilon = \widetilde{\Upsilon}^{-1}$.

%Consider
{Applying the Finsler's Lemma~\cite{Pipeleers:2009aa}, consider:}
\begin{align}\label{eq:proof-lemma4-3}
    \mathcal{B} = \begin{bmatrix}
        I  & 0 \\ H K & - H \\ 0 & I
    \end{bmatrix}.
\end{align}
By multiplying~\eqref{eq:proof-lemma4-2} on the left by $\mathcal{B}^\top$ on the left and its transpose on the right, we obtain
\begin{align}\label{eq:proof-lemma4-4}
    \begin{bmatrix}
        PHK + K^\top H P + 2 \eta P & L^\top \Upsilon - PH \\
        \Upsilon L - HP & -2 \Upsilon
    \end{bmatrix} < 0.
\end{align}
By multiplying~\eqref{eq:proof-lemma4-4} on the left by $[\hat{G}_{\mathrm{av}}^\top(t) \; \psi^\top(K \hat{G}_{\mathrm{av}}(t))]$ and its transpose on the right, yields
%\begin{small}
\begin{align}
     &\hat{G}_{\mathrm{av}}^\top(t) \left(PHK + K^\top H P\right)\hat{G}_{\mathrm{av}}(t) - 2 \hat{G}_{\mathrm{av}}^\top P H \psi(K \hat{G}_{\mathrm{av}}(t)) \nonumber \\
     &- 2 \psi^\top(K \hat{G}_{\mathrm{av}}) \Upsilon \left(\psi(K \hat{G}_{\mathrm{av}}) \!-\! L \hat{G}_{\mathrm{av}}\right) \nonumber\\
     &+ 2 \eta \hat{G}_{\mathrm{av}}^\top(t) Q \hat{G}_{\mathrm{av}}(t) \!<\! 0.
\end{align}
%\end{small}
Now, by multiplying the inequalities in~\eqref{eq:saturation_thm} on the left by $\mathrm{diag}(X^{-\top}, 1)$ and its tranponse on the right, we obtain
\begin{align}\label{eq:proof-lemma4-5}
     \begin{bmatrix}
            P &  K_{(\ell)}^\top - L_{(\ell)}^\top \\
            K_{(\ell)} - L_{(\ell)} & \overline{u}_\ell^2
    \end{bmatrix} \geq 0, \quad \ell = 1,\ldots,n.
\end{align}
From the Schur complement lemma, we have that~\eqref{eq:proof-lemma4-5} implies
\begin{align}\label{eq:proof-lemma4-6}
    P - \frac{1}{\overline{u}_\ell^2} (K_{(\ell)} - L_{(\ell)})^\top (K_{(\ell)} - L_{(\ell)}) \geq 0, \; \ell = 1,\ldots,n.
\end{align}
By multiplying~\eqref{eq:proof-lemma4-6} {on the left by $\hat{G}_{\mathrm{av}}^\top$} and its transpose on the right, we get
\begin{align}
     V(\hat{G}_{\mathrm{av}}) \geq \frac{|(K-L)_{(\ell)} \hat{G}_{\mathrm{av}}|^2}{\overline{u}_{(\ell)}^2}, \; \ell = 1,\ldots,n.
\end{align}
Then, provided that $\hat{G}_{\mathrm{av}} \in \mathcal{E}$, we ensure that $\hat{G}_{\mathrm{av}} \in \mathcal{G}$, that is, $\mathcal{E} \subset \mathcal{G}$, and the conditions of Lemma~\ref{lem:saturacao} are satisfied. It allows us to obtain from~\eqref{eq:proof-lemma4-4} and~\eqref{eq:setor_psi} that
\begin{align}\label{eq:exp-ineq}
    \dot{V}(\hat{G}_{\mathrm{av}}(t)) \leq - 2 \eta V(\hat{G}_{\mathrm{av}}(t)) < 0, \, \forall \hat{G}_{\mathrm{av}}(t) \neq 0,
\end{align}
where $V(\hat{G}_{\mathrm{av}})$, defined in~\eqref{eq:lyapunov_Gav},
is a Lyapunov function that ensures the exponential stability of the origin of the average system.
From the Comparison Lemma, it follows from~\eqref{eq:exp-ineq} that
\begin{align}
    {V}(\hat{G}_{\mathrm{av}}(t)) \leq e^{-2 \eta t} {V}(\hat{G}_{\mathrm{av}}(0)).
\end{align}
Furthermore, as 
\begin{align}
    \lambda_{\min}(P) \|\hat{G}_{\mathrm{av}}\|^2 \leq 
    V(\hat{G}_{\mathrm{av}}) \leq \lambda_{\max}(P) \|\hat{G}_{\mathrm{av}}\|^2,
\end{align}
we can obtain
\begin{align}
    \|\hat{G}_{\mathrm{av}}(t)\| \leq \kappa_g e^{-\eta t}\|\hat{G}_{\mathrm{av}}(0)\|
\end{align}
where $\kappa_g = \sqrt{\lambda_{\max}(P) / \lambda_{\min}(P)}$.
This concludes the proof.~\hfill$\blacksquare$
\end{pf}

\subsubsection{Practical Exponential Stability via Averaging Theorem}

{Lemma~\ref{lem:Gav_stability} gives a condition to design the control gain ensuring regional exponential stability of the average closed-loop system~\eqref{eq:dinamica_Gav_politopo}. Next, we present the main result, showing local exponential convergence to a neighborhood of the extremum via Averaging Theory (see Appendix~\ref{sec:appendix}).}

\begin{theorem}
    Consider the ESC system in Fig.~\ref{fig:diag_ES_sat} with locally quadratic nonlinear map (\ref{eq:map_static}),  Assumptions~\ref{assump:probing_w} and \ref{assump:H*} as well as   
    %Consider
    the corresponding average closed-loop dynamics 
    governing the gradient estimate subject to saturation in~\eqref{eq:average_system1}.
    If the conditions of Lemma~\ref{lem:Gav_stability} are satisfied, 
    then, for $\omega > 0$ sufficiently large in (\ref{eq:relacao_sinais_omega_pert_input}), 
    %the equilibrium $\hat{G}_{\mathrm{av}} = 0$ is exponentially stable and %$\tilde{\theta}_{\mathrm{av}}(t)$ converges exponentially to zero. In %particular,
    there exist constants $\bar{\kappa}_\theta, \bar{\kappa}_y, \eta > 0$,  such that: 
    \begin{align} \label{eq:theta_desigualdades}
        \|\theta(t) - \theta^\ast\| \leq \bar{\kappa}_\theta e^{-\eta t} + \mathcal{O}\left(a + \frac{1}{\omega}\right), \\
        \label{eq:y_desigualdades}
        |y(t) - Q^\ast| \leq  \bar{\kappa}_y e^{-\eta t} + \mathcal{O}\left(a^2 + \frac{1}{\omega^2}\right),
    \end{align}
    where $a = \sqrt{\sum_{i=1}^n a_i^2}$, with $a_i$ defined in  (\ref{eq:vetor_S_input_contrained})--(\ref{eq:vetor_M_input_constrained}), and $\bar{\kappa}_\theta$ and $\bar{\kappa}_y$ are constants which depend on the initial condition $\theta(0)$. 
\end{theorem}
\begin{pf}
    From equations~\eqref{eq:grad_simp_2} and \eqref{eq:delta_input_constrained}, and recalling that $\Delta_{\mathrm{av}}(t)=0$, it can be obtained that
    \begin{align}
        \hat{G}_{\mathrm{av}}(t) = H \tilde{\theta}_{\mathrm{av}}(t),
    \end{align}
    since 
    %the quadratic term $\frac{1}{2}M(t)\tilde{\theta}^\top H\tilde{\theta}%(t)$ can be neglected in a local analysis, and
    the other terms also have zero mean.
    
    Rewriting the Lyapunov function in~\eqref{eq:lyapunov_Gav} as
    \begin{align}
        V(\tilde{\theta}_{\mathrm{av}}) = \tilde{\theta}_{\mathrm{av}}^\top \overline{P} \tilde{\theta}_{\mathrm{av}},
    \end{align}
    where $\overline{P} = H P H$ is a symmetric positive definite matrix, provided that {$P$ is symmetric and positive definite and $H$ is symmetric}. Thus, it is possible to find
    \begin{align} \label{palmeirasevasco}
        \|\tilde{\theta}_{\mathrm{av}}(t)\| \leq \kappa_\theta e^{-\eta t} \|\tilde{\theta}_{\mathrm{av}}(0)\|,
    \end{align}
    where $\kappa_\theta = \sqrt{\lambda_{\max}(\overline{P}) / \lambda_{\min}(\overline{P})}$. Since the closed-loop average system is exponentially stable according to (\ref{palmeirasevasco}), by applying averaging theorem \cite{plotnikov1979averaging} (see Appendix~\ref{sec:appendix}, with $\varepsilon:=1/\omega$), it follows that:   
%As the differential equation has discontinuity on the right side, due to the presence of the saturating function, \eqref{eq:esc_temp_mudanca_Grad_function} is $T$-periodic and Lipschitz continuous, it follows from~\cite{plotnikov1979averaging} that 
    % $\tilde{\theta}_{\mathrm{av}}(\tau)$ is asymptotically stable and, consequently, it can be guaranteed that
    \begin{align}
        \|\tilde{\theta}(t) - \tilde{\theta}_{\mathrm{av}}(t)\| \leq \mathcal{O}\left(\frac{1}{\omega}\right).
    \end{align}
    Applying the triangle inequality, it is guaranteed that
    \begin{align} \label{meucudeouro}
        \| \tilde{\theta}(t) \| 
        &\leq \kappa_\theta e^{-\eta t} \| \tilde{\theta}_{\text{av}}(0) \| + \mathcal{O}\left(\frac{1}{\omega}\right).
    \end{align}
    From the averaging theorem~\cite{plotnikov1979averaging}, it can also be concluded that
    \begin{align}
        \| \hat{G}(t) - \hat{G}_{\text{av}}(t) \| \leq \mathcal{O}\left(\frac{1}{\omega}\right).
    \end{align}
    Similarly, we can apply the triangle inequality to obtain
    \begin{align}
        \| \hat{G}(t) \| 
        % &\leq \| \hat{G}_{\text{av}}(t) \| + \mathcal{O}\left(\frac{1}{\omega}\right)  \\
        &\leq \kappa_{g} e^{-\eta t} \| \hat{G}_{\text{av}}(0) \| + \mathcal{O}\left(\frac{1}{\omega}\right).
    \end{align}
    From~\eqref{relacao_theta_til&theta_medio} and (\ref{meucudeouro}), the following relation can be obtained:
    \begin{align} \label{eq:demonstracao_estabilidade_theta}
        \| \theta(t) - \theta^* \| \leq \kappa_\theta e^{-\eta t} \|\theta(0)- \theta^*\| + \mathcal{O}\left(a + \frac{1}{\omega}\right),
    \end{align}
resulting in (\ref{eq:theta_desigualdades}), with $\bar{\kappa}_\theta=\kappa_\theta\|\theta(0)- \theta^*\|$. 

Let the output error be 
\begin{align}
    \tilde{y}(t) := y(t) - Q^*,~~~~y(t) = Q(\theta(t)).
\end{align}
By computing its norm and using the Cauchy–Schwarz inequality, one gets
\begin{align}
    |\tilde{y}(t)| & = |y(t) - Q^*| = |(\theta(t) - \theta^*)^\top H(\theta(t) - \theta^*)| \\ &\leq \|H\| \|\theta(t) - \theta^*\|^2.
\end{align}
From~\eqref{eq:demonstracao_estabilidade_theta}, it is still possible to obtain
\begin{align}
    |\tilde{y}(t)| \!\leq\! \|H\| \kappa_\theta^2 e^{-2\eta t} \|\theta(0)\!-\!\theta^*\|^2 \!+\! \mathcal{O}\left(a^2 \!+\! \frac{2a}{\omega} \!+\! \frac{1}{\omega^2}\right).
\end{align}
Since $e^{-\eta t} \geq e^{-2 \eta t}$ and $a^2+\frac{1}{\omega^2} \geq \frac{2a}{\omega}$, for $\omega > 0$ and $a>0$, by the Young's inequality, one obtains
\begin{align}
    |y(t) - Q^*| \leq \bar{\kappa}_y e^{-\eta t} + \mathcal{O}\left(a^2 + \frac{1}{\omega^2}\right),
\end{align}
where 
\begin{align*}
    \bar{\kappa}_y = \kappa_\theta^2 \|H\| \|\theta(0)- \theta^*\|^2, 
\end{align*}
resulting in inequality \eqref{eq:y_desigualdades}. This concludes the proof. ~\hfill$\blacksquare$
\end{pf}

%---------------------------------------------------%

%%%%% Artigo com a Saturação no Gradiente %%%%%

\section{Numerical Results}

The effectiveness of the proposed approaches is illustrated via two numerical examples. 
The method for ESC under input saturation from Section~\ref{sec:preliminaries} is validated in Example~1, 
while the ESC \textcolor{black}{under gradient saturation} from Section 3 is addressed in Example~2.

\subsection{Example~1: ESC under Input Saturation}

Consider the ESC system under input saturation with a nonlinear map \eqref{eq:map_static_input_constrained} 
with an unknown Hessian matrix taking values in the polytopic domain given by the following vertices
\begin{align}
    H_1 = (1-\overline{\delta}) H_0, \quad  H_2 = (1+\overline{\delta}) H_0,
\end{align}
where $\overline{\delta} > 0$ is a parameter used to construct the vertices of the polytopic domain and $H_0$ is the Hessian matrix used in~\cite{ghaffari2012multivariable}:
\begin{align}
    H_0 = \begin{bmatrix}
        100 & 30 \\
        30 & 20
    \end{bmatrix} > 0.
\end{align}
In addition, for the simulations, it is assumed that unknown parameters are $Q^* = 10$ and $\theta^* = [2 \; 4]^\top$. 
Note that the unknown parameters $Q^*$ and $\theta^*$ are not used to design the gains of the AW controller~\eqref{eq:control_law_1}.
For illustration purposes, we assume that $\overline{\delta} = 0.1$.
Then, the gains of the AW controller~\eqref{eq:control_law_1} were designed by solving the conditions in Lemma~\ref{thm:1_input_constrained} with 
the decay rate $\eta = 1$ and saturation levels $\overline{\theta}_1 = \overline{\theta}_2 = 5$. 
The resulting control gains are the following:
$$K = \begin{bmatrix}
     %  -0.0281  &  0.0405 \\
    %0.0442 &  -0.1463
       -0.0270 &   0.0361 \\
    0.0456  & -0.1492
\end{bmatrix}, \, K_{\mathrm{aw}} = \begin{bmatrix}
    %1.5534  &  0.0329 \\
   %-0.0332  &  1.5533
       2.2794  & 0.0824 \\
       -0.0865  &  2.2804
\end{bmatrix}.$$ 

\begin{figure}[ht!]
     \centering
%------------Figura sinal de theta------------
    \begin{subfigure}[b]{\columnwidth}
        \centering
        \includegraphics[width=0.8\textwidth]{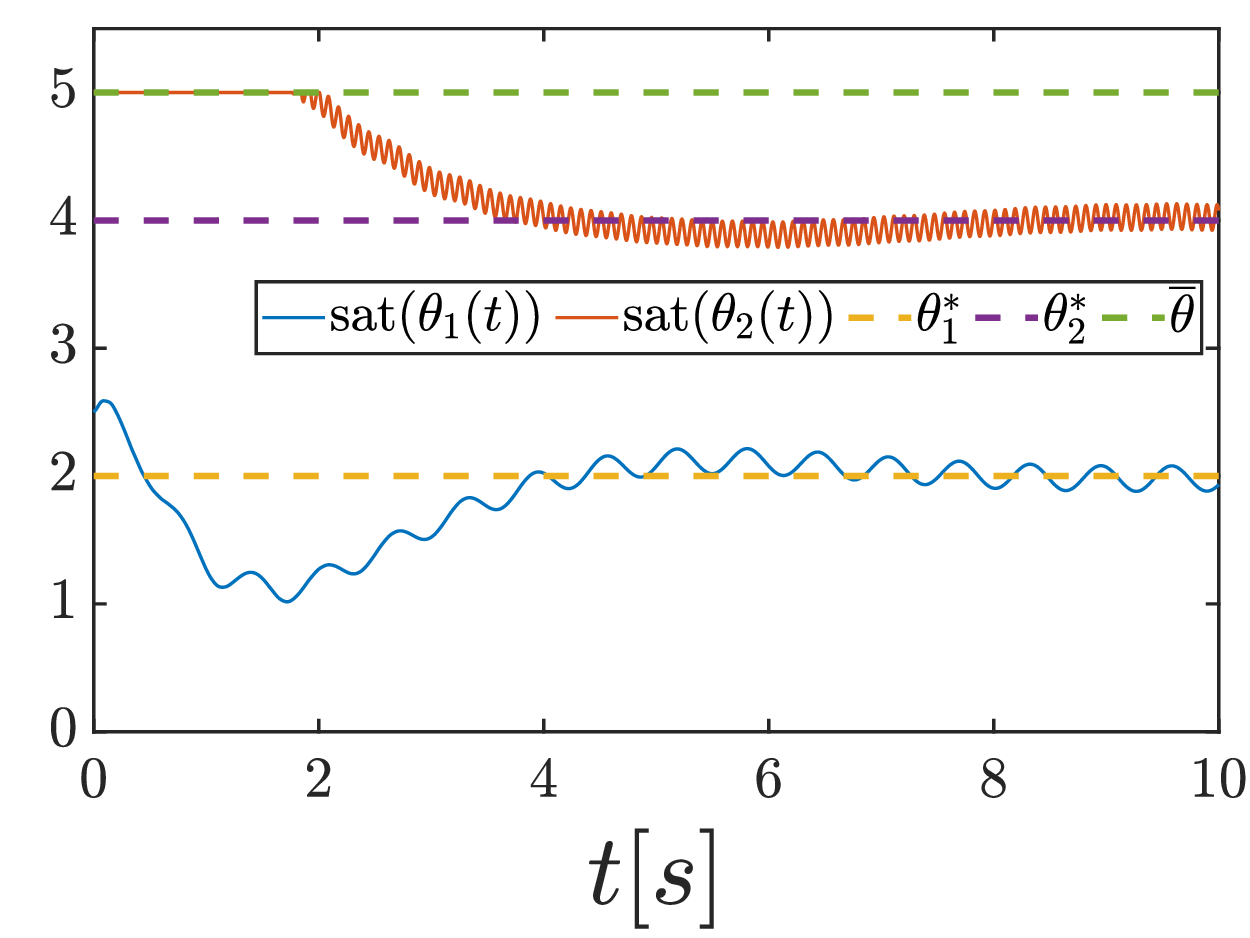}
        \caption{Saturated input of the map -- $\mathrm{sat}(\theta(t))$.}
    \end{subfigure}
%------------Figura sinal da saída------------
    \begin{subfigure}[b]{\columnwidth}
        \centering
        \includegraphics[width=0.8\textwidth]{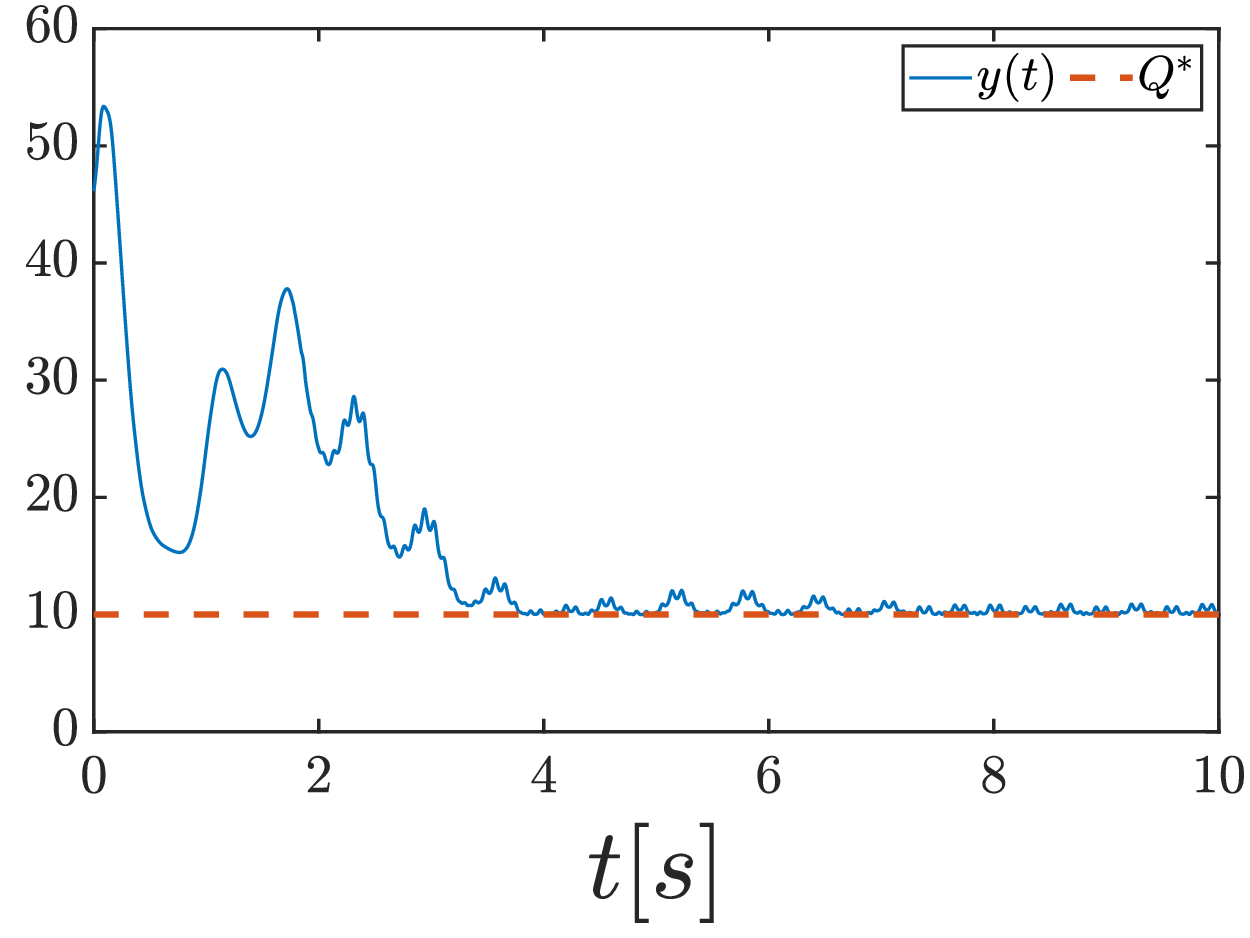}
        \caption{Output of the map -- $y(t)$.}
    \end{subfigure}
\caption{Trajectories of the closed-loop ESC system under input saturation with anti-windup controller~\eqref{eq:control_law_1} designed according to Lemma~\ref{thm:1_input_constrained} -- Example~1.}
\label{fig:example1}
\end{figure}

For the simulations, the dither signals $S(t)$ and $M(t)$ given in~\eqref{eq:vetor_S_input_contrained} and~\eqref{eq:vetor_M_input_constrained}, respectively, are selected with frequencies $\omega_1 = 10$ rad/s and $\omega_2 = 70~\mathrm{rad/s}$, and amplitudes $a_1=a_2=0.1$. Besides that, the simulations were performed considering the initial condition $\theta(0) =  [2.5 \; 6]^\top$ and $\alpha = [0.6822 \;0.3178]^\top$.
The trajectories of the closed-loop ESC under input saturation with the designed AW controller~\eqref{eq:control_law_1} 
are shown in Fig.~\ref{fig:example1}. In Fig.~\ref{fig:example1}(a), it is possible to notice the convergence of the inputs to the neighborhood of the optimum point $\theta^*$, even in the presence of saturation. 
Moreover, the convergence of the output to the neighborhood of $Q^*$ is shown in~Fig.~\ref{fig:example1}(b).
This clearly illustrates the theoretical findings established in Theorem~\ref{thm:main_result_1}.

\begin{figure}[ht!]
     \centering
%------------Figura sinal de theta------------
    \begin{subfigure}[b]{\columnwidth}
        \centering
        \includegraphics[width=0.8\textwidth]{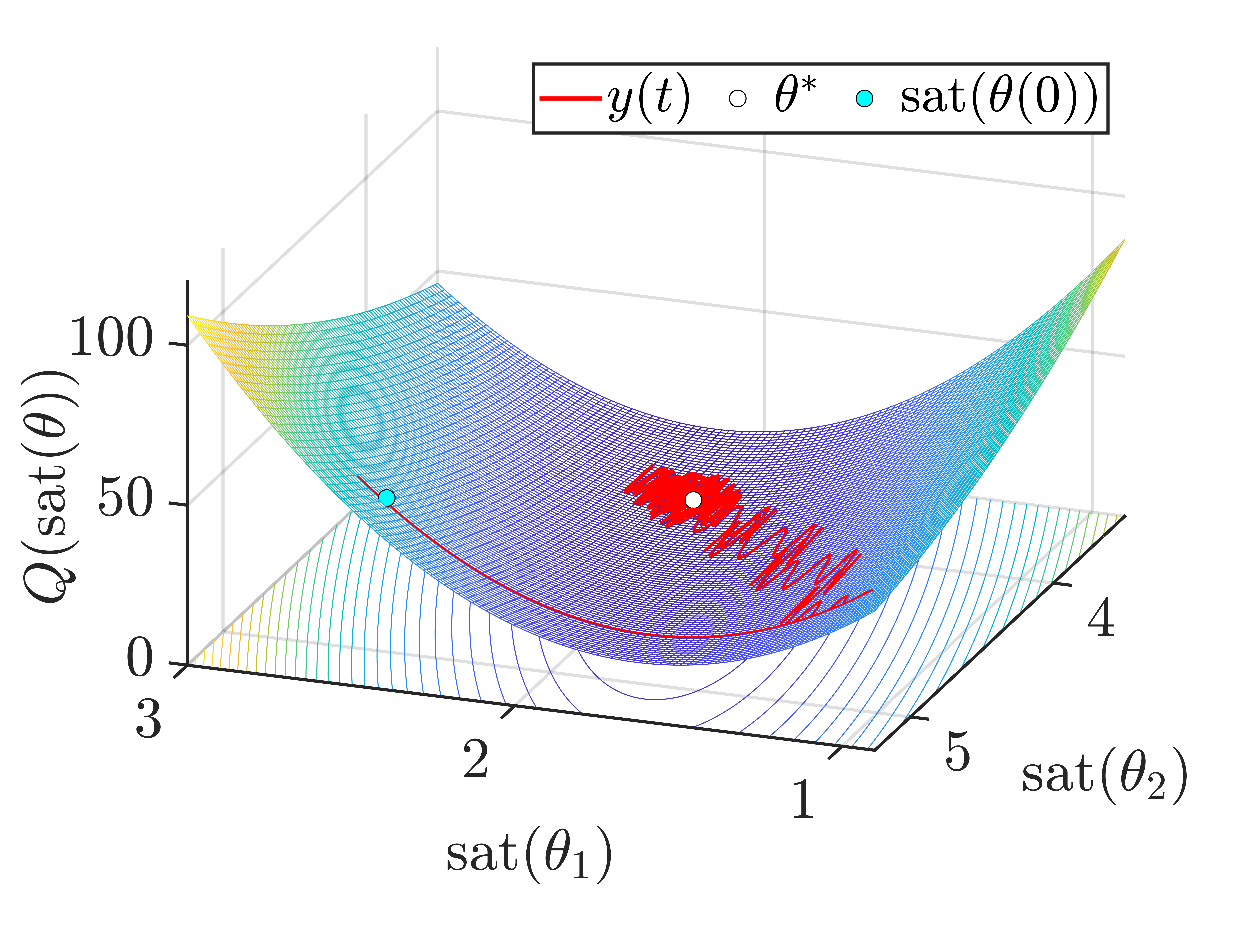}
        \caption{Trajectory $y(t)$ along the surface of the quadratic map with input saturation $Q(\mathrm{sat}(\theta))$ in~\eqref{eq:map_static_input_constrained}.}
    \end{subfigure}
% \subfloat[\label{fig:surface_side_view}Trajectory $y$ over the paraboloid surface]{
%     \includegraphics[width=0.4\textwidth]{Graph_3D_antiwindup/trajetoria_y_3D_lateral.eps}}
%     \hfill
% \vskip 3mm
%------------Figura sinal da saída------------
    \begin{subfigure}[b]{\columnwidth}
        \centering
        \includegraphics[width=0.8\textwidth]{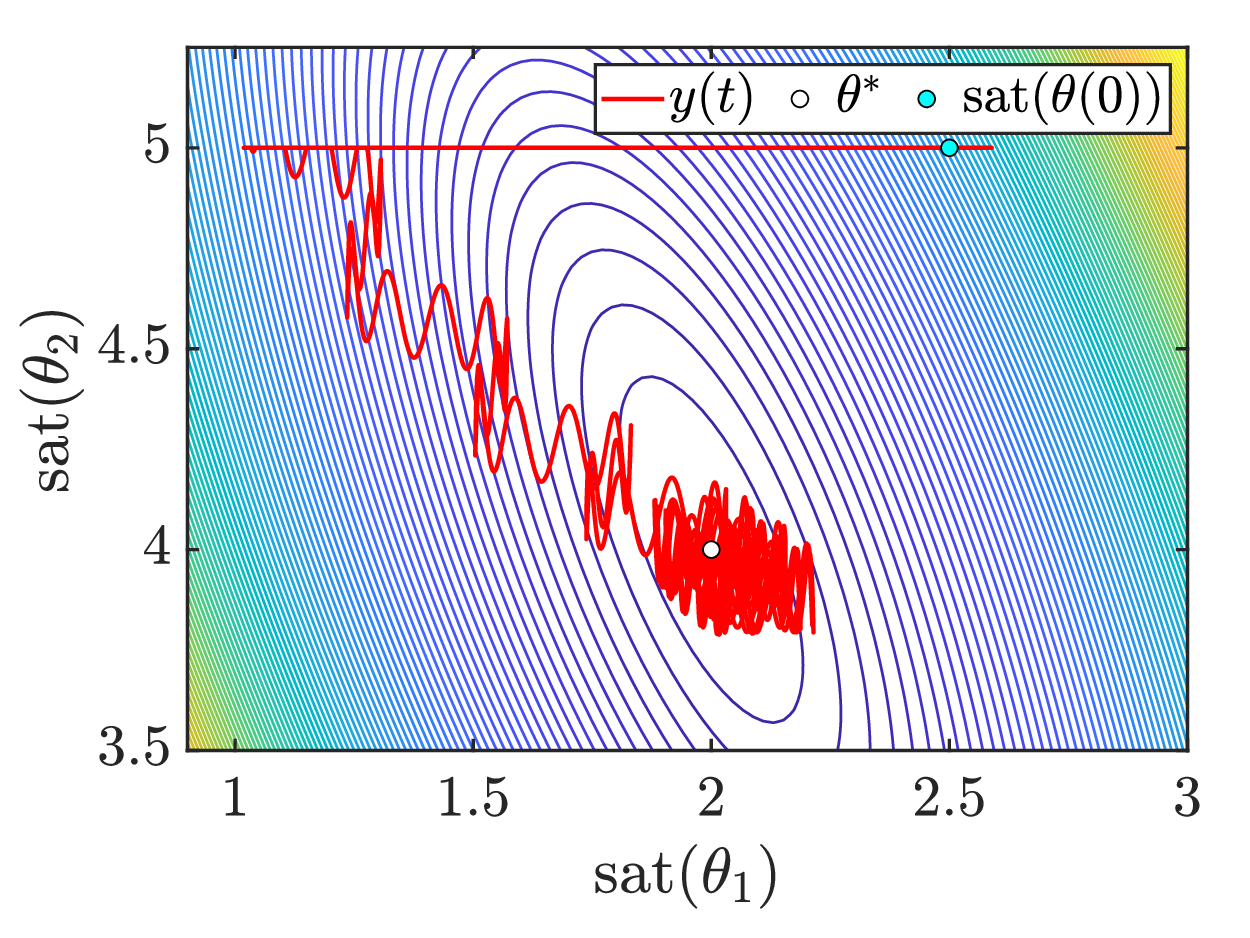}
        \caption{Trajectory $\mathrm{sat}(\theta(t))$ and level sets of the quadratic map with input saturation $Q(\mathrm{sat}(\theta))$ in~\eqref{eq:map_static_input_constrained}.}
    \end{subfigure}
    % \subfloat[\label{fig:surface_top_view}Evolution of the saturated inputs along the level surfaces.]{
    % \includegraphics[width=0.4\textwidth]{Graph_3D_antiwindup/trajetoria_y_3D_superior.eps}}
    % \hfill  
\caption{Trajectory of the output $y(t)$ of the closed-loop ESC system under input saturation with the anti-windup controller~\eqref{eq:control_law_1} designed according to Lemma~\ref{thm:1_input_constrained} -- Example~1.}
\label{fig:surface_y}
\end{figure}

Furthermore, Fig.~\ref{fig:surface_y}(a) depicts the evolution of $y(t)$ along the surface of the quadratic map with input saturation $Q(\mathrm{sat}(\theta))$ in~\eqref{eq:map_static_input_constrained}. The trajectory of $\mathrm{sat}(\theta(t))$ together with several level sets of $Q(\mathrm{sat}(\theta))$ are shown in~Fig.~\ref{fig:surface_y}(b).
%%%%%
%%%%%
%In the Figure~\ref{fig:comparacao_ganho_LMI_arbitrario_Diag_input}, the extremum-seeking control system with saturation was simulated with the designed controller using the Theorem conditions~\ref{thm_input_constrained:3}, developed in this work, and with a negative diagonal gain given by $K = -0.025 I_2$ and without the Anti-windup compensator.

\begin{figure}[ht!]
     \centering
%------------Figura sinal de theta------------
    \begin{subfigure}[b]{\columnwidth}
        \centering
        \includegraphics[width=0.8\textwidth]{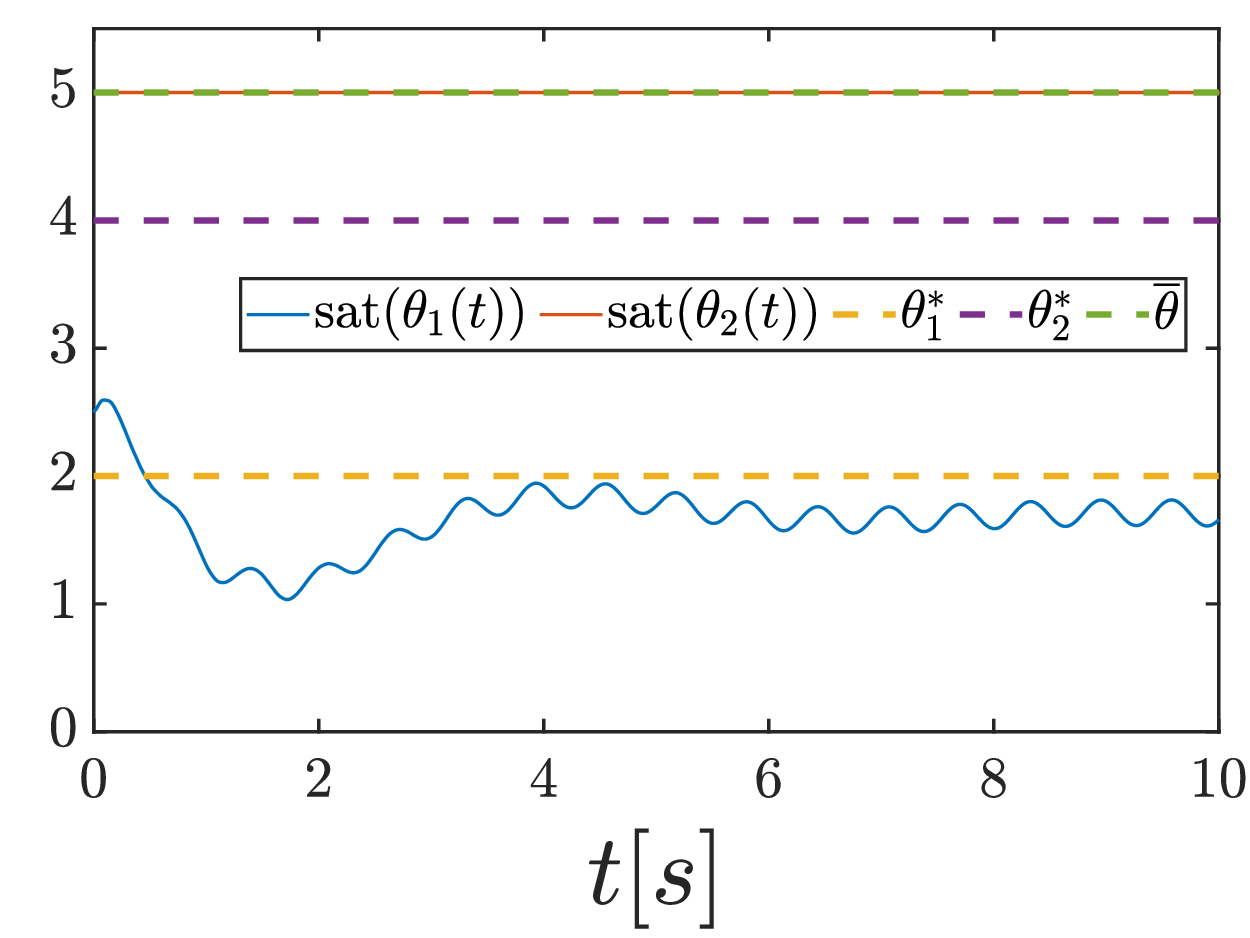}
        \caption{Saturated input of the map -- $\mathrm{sat}(\theta(t))$.}
    \end{subfigure}
%------------Figura sinal da saída------------
    \begin{subfigure}[b]{\columnwidth}
        \centering
        \includegraphics[width=0.8\textwidth]{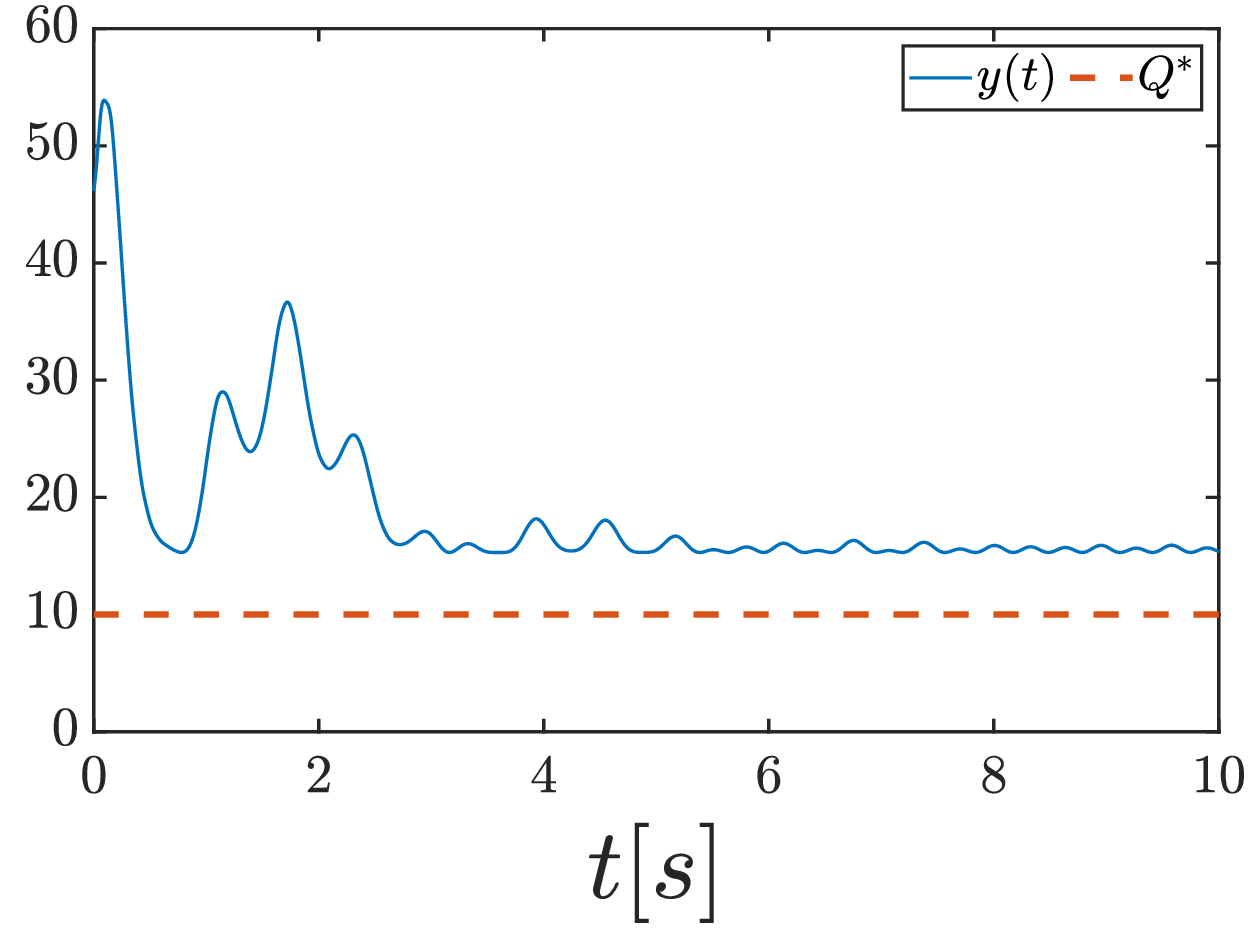}
        \caption{Output of the map -- $y(t)$.}
    \end{subfigure}
\caption{Trajectories of the closed-loop ESC system under input saturation without anti-windup compensation -- Example~1.}
\label{fig:example1_withoutAW}
\end{figure}

Consider the control law~\eqref{eq:control_law_1} without the AW compensation term, that is, $K_{\mathrm{aw}} = 0$. The same control gain $K$ is considered. The closed-loop simulation for this case is shown in Fig.~\ref{fig:example1_withoutAW}.
{In Fig.~\ref{fig:example1_withoutAW}, the ESC system does not converge to the extremum, highlighting the advantages of the proposed AW-based approach shown in Fig~\ref{fig:example1}, which guarantees convergence.}

\subsection{Example~2: ESC under \textcolor{black}{Gradient Saturation}}
\begin{figure}[ht!]
     \centering
%------------Figura sinal de controle-----------
    \begin{subfigure}[b]{\columnwidth}
        \centering
        \includegraphics[width=0.8\textwidth]{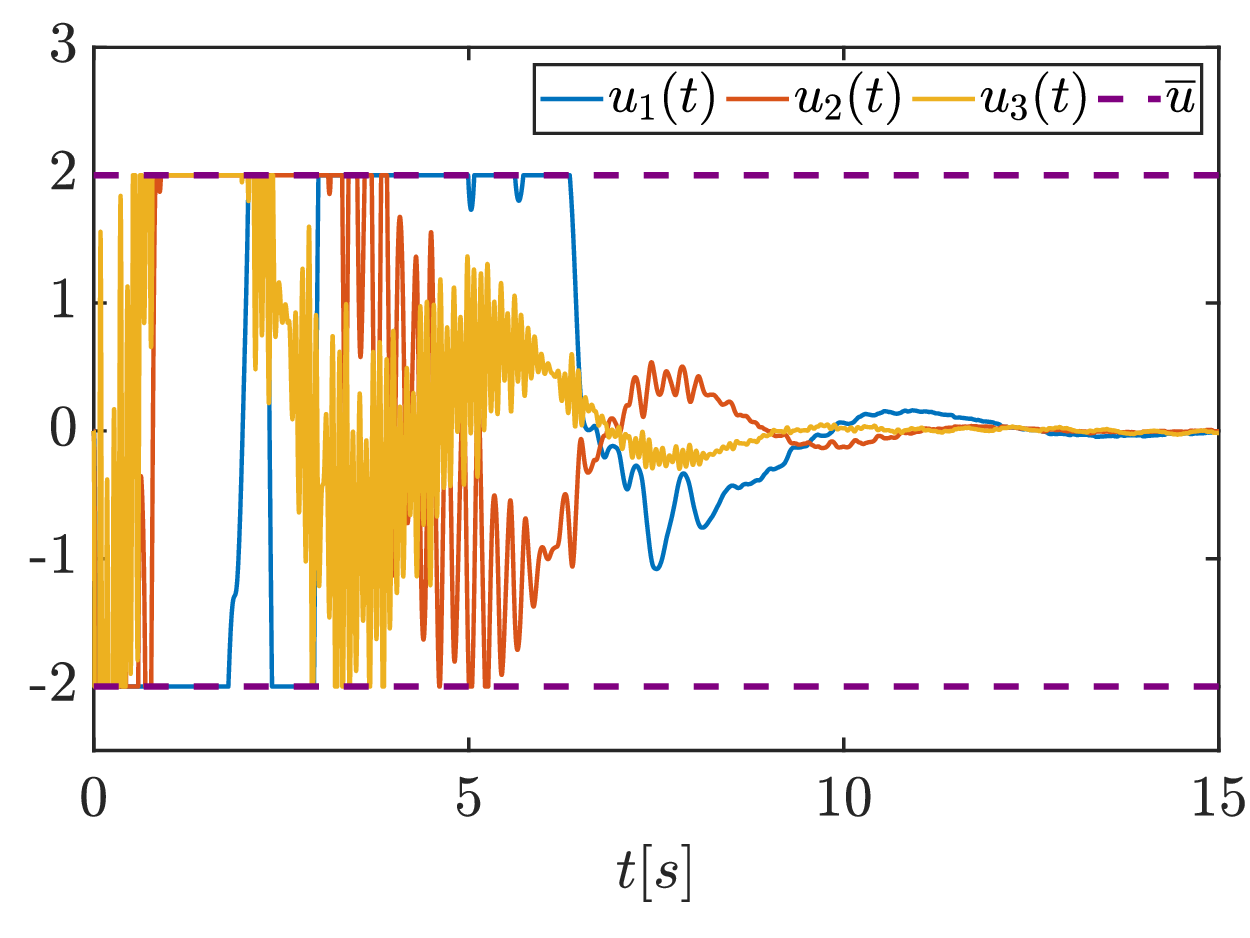}
        \caption{\textcolor{black}{Control signal with gradient saturation -- $u(t) = \mathrm{sat}(K \hat{G}(t))$.}}
    \end{subfigure}
% \subfloat[\label{fig:sinal_controle_LMI_dim3}$\mathrm{sat}(u(t))$ -- 3 dimensions]{
%     \includegraphics[width=0.4\textwidth]{Simulacoes_sat_gradient_3dim/sinal de sat(u(t)).eps}}
%     % \hfill
% % \vskip 3mm
% %------------Figura sinal de theta------------
\begin{subfigure}[b]{\columnwidth}
        \centering
        \includegraphics[width=0.8\textwidth]{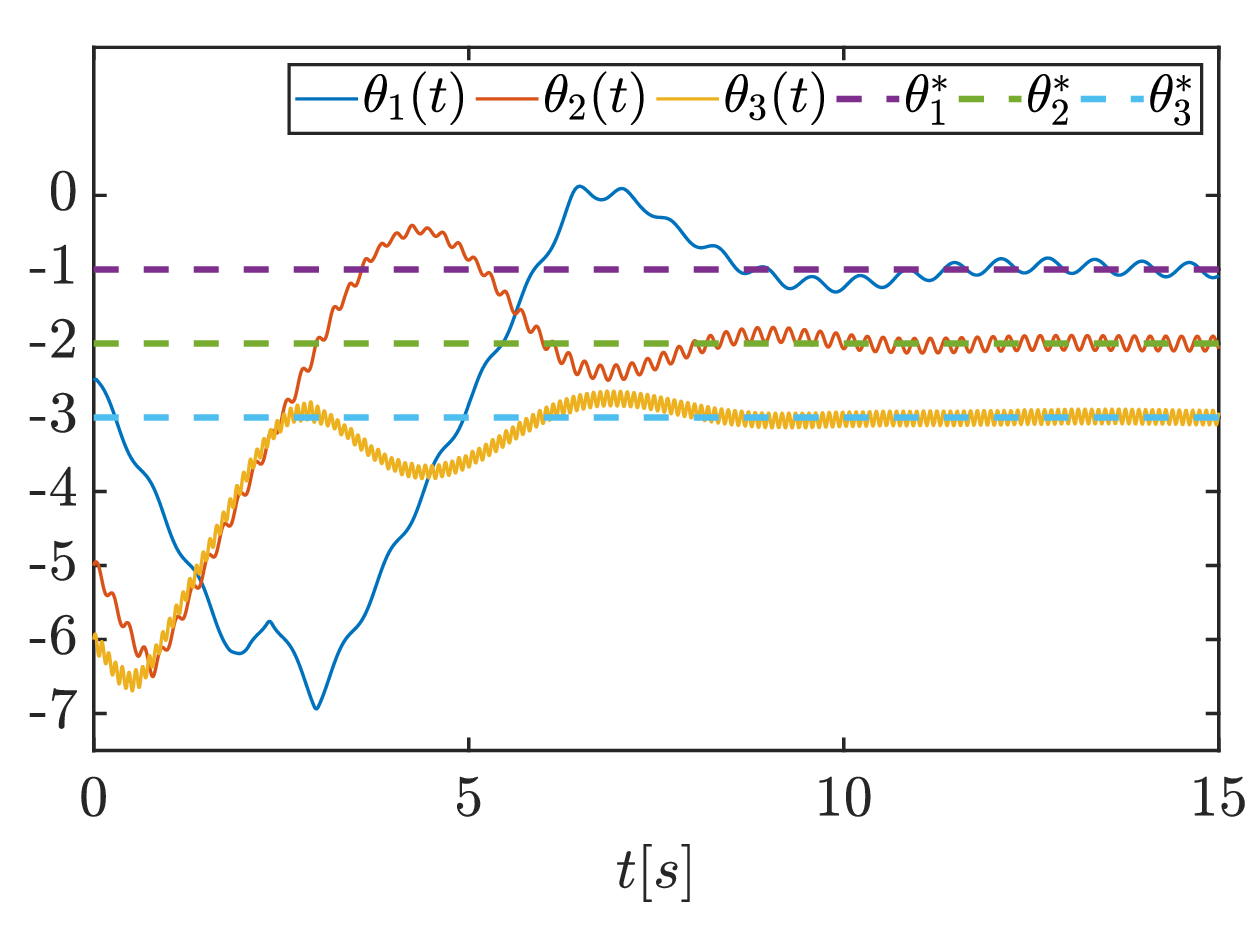}
        \caption{Input vector of the map -- $\theta(t)$.}
\end{subfigure}
%     \subfloat[\label{fig:sinal_theta_LMI_dim3}$\theta(t)$ -- 3 dimensions]{
%     \includegraphics[width=0.4\textwidth]{Simulacoes_sat_gradient_3dim/sinal de theta(t).eps}}
%     % \hfill
% % \vskip 3mm
% %------------Figura sinal da saída------------
\begin{subfigure}[b]{\columnwidth}
        \centering
        \includegraphics[width=0.8\textwidth]{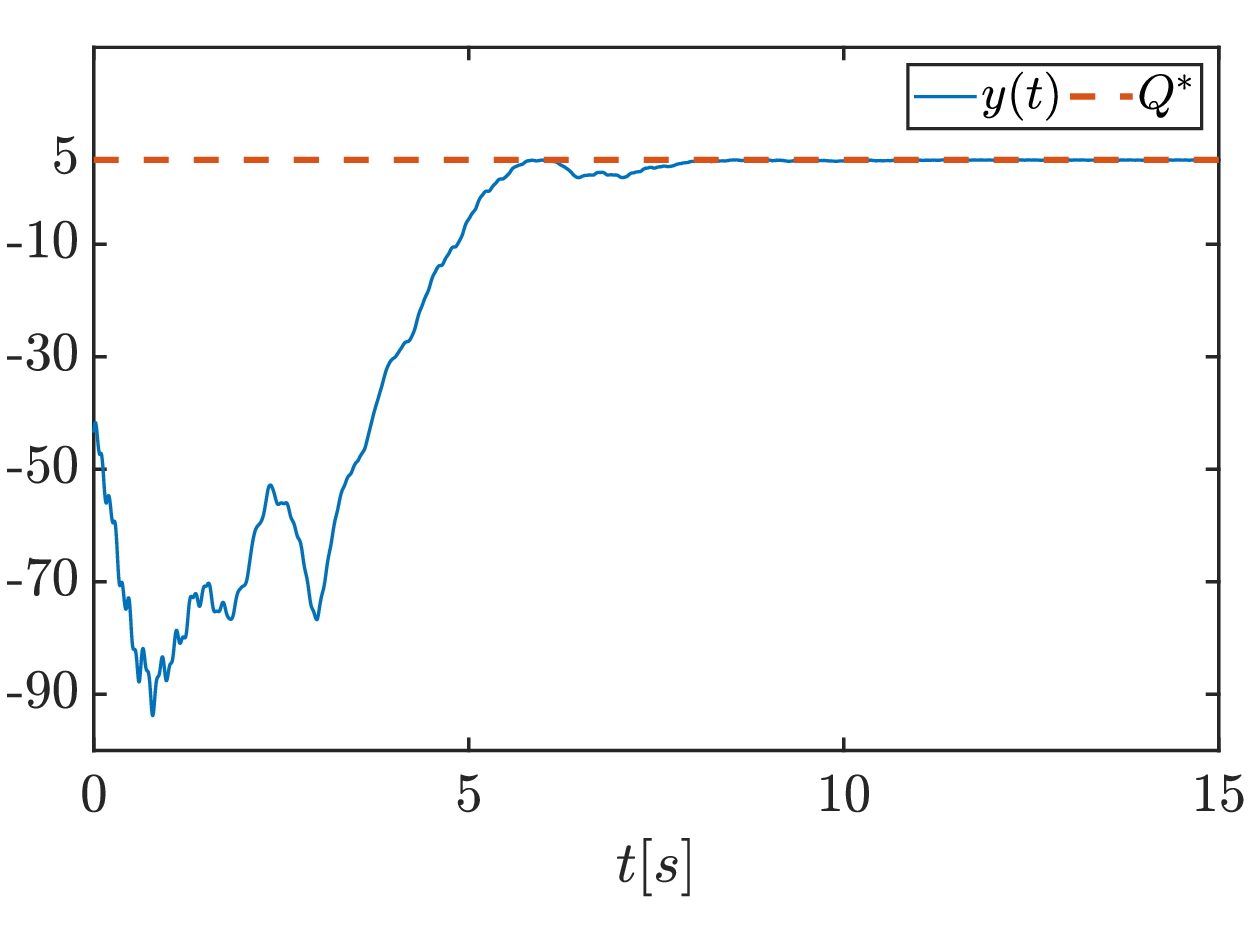}
        \caption{Output of the map -- $y(t)$.}
\end{subfigure}
%     \subfloat[\label{fig:sinal_saida_LMI_dim3}$y(t)$ -- 3 dimensions]{
%     \includegraphics[width=0.4\textwidth]{Simulacoes_sat_gradient_3dim/sinal de y.eps}}
%     % \hfill    
\caption{Trajectories of the closed-loop ESC system \textcolor{black}{under gradient saturation} with controller~\eqref{eq:control_law_1_cosmeedamiao} designed according to Lemma~\ref{lem:Gav_stability} -- Example~2.}
\label{fig:example2}
\end{figure}

Consider the ESC system under gradient saturation discussed in Section~\ref{sec:3}.
We consider a nonlinear map in~\eqref{eq:map_static} with three inputs and unknown parameters $Q^* = 5$ and 
$\theta^* = [-1 \; -2 \; -3]^\top$. 
In this case, we consider the uncertain Hessian matrix taken in a polytopic domain with four {negative definite vertices randomly generated} $H \in \mathrm{co}\{H_1,H_2,H_3,H_4\}$, where
\begin{align*}
    H_1 = 
    \begin{bmatrix}
        -6.7828  &  0.8480 &  -1.3462 \\
        0.8480 &  -6.0017 &  -0.7825 \\
       -1.3462  & -0.7825 &  -3.2421
    \end{bmatrix}, 
    \end{align*}
    \begin{align*}
    H_2 = 
    \begin{bmatrix}
          -3.9159 &  -0.8122  &  1.4150 \\
       -0.8122 &  -5.7484 &  -0.0047 \\
        1.4150  & -0.0047 &  -4.6956
    \end{bmatrix}, \\
    H_3 = 
    \begin{bmatrix}
        -3.9141 &  -0.3951 &   0.5802 \\
       -0.3951 &  -3.6059 &   1.0325 \\
        0.5802 &   1.0325 &  -4.0962
    \end{bmatrix}, \\
    H_4 = 
    \begin{bmatrix}
        -6.1443  &  0.0911  & -0.7984 \\
        0.0911 &  -5.9879 &  -2.3066 \\
       -0.7984 &  -2.3066  & -3.9025
    \end{bmatrix}.
\end{align*}
The control gain is designed by solving the conditions in Lemma~\ref{lem:Gav_stability} with
$\epsilon = 0.5$, decay rate $\eta=1$, and the saturation levels are $\overline{u}_1 = \overline{u}_2 = \overline{u}_3 = 2$.
The designed control gain is
$$K = \begin{bmatrix}
    0.5009 &  -0.0094 &  -0.0018 \\
   -0.0104  &  0.5312 &  -0.0881 \\
    0.0006 &  -0.0856  &  0.7352
\end{bmatrix}.$$
To perform the closed-loop simulation, we consider that the dither frequencies are $\omega_1 = 10$ rad/s, $\omega_2 = 30~\mathrm{rad/s}$ and $\omega_3 = 70~\mathrm{rad/s}$, their amplitudes are $a_1=a_2=a_3=0.1$, and the initial condition is $\theta(0) =  [2.5 \; 5\;6]^\top$.
The results obtained with the closed-loop simulation are shown in Fig.~\ref{fig:example2}. 
Particularly, Fig.~\ref{fig:example2}(a) depicts the trajectory of the \textcolor{black}{$u(t)$}.
It is possible to notice that the signal $u(t)$ converges exponentially to zero, indicating the convergence of
the gradient estimate to zero, even in the presence of saturation. As a result, the input of the quadratic map
converges to the neighborhood of the unknown point $\theta^*$, as shown in Fig.~\ref{fig:example2}(b), and the output $y(t)$
converges to the the neighborhood of the optimum point $Q^* = 5$, as shown Fig.~\ref{fig:example2}(c).

\section{Conclusion}

This paper has addressed the challenging problem of multivariable extremum seeking control in the presence of both actuator saturation and gradient saturation. By employing a sector representation, we established stability analysis conditions for the average system under these saturation effects, thereby extending the applicability of extremum seeking to more realistic scenarios where input and gradient constraints cannot be ignored. To rigorously justify the stability claims, the averaging theorem for non-differentiable Lipschitz systems was invoked, ensuring that the trajectories of the closed-loop system converge to a neighborhood of the unknown optimal point, even under limited actuation and bounded gradient information. Moreover, by assuming an uncertain polytopic representation of the Hessian matrix, constructive and \textcolor{black}{verifiable LMI conditions} were derived for designing stabilizing controllers, providing a systematic framework that can be applied to a broad class of nonlinear optimization problems. Numerical simulations further illustrated the practicality and effectiveness of the proposed feedback controllers by demonstrating the convergence of the system to the extremum point, confirming the robustness of the design against uncertainties and saturation effects. Overall, the results presented in this work offer both theoretical insights and practical tools for extremum seeking control in constrained multivariable settings.

Future investigation lies in the expansion of the proposed design and analysis, taking into account saturation constraints for different control problems with unknown control direction, and pursuing infinite-dimensional multi-agent optimization via
Nash equilibrium seeking, as considered in references \cite{IJACSP:2017,JOTA:2021}.
{Other directions include the proposal of AW compensation of the extremum seeking with bounded update rates with potential applications in aeronautical systems~\cite{biannic2009optimization} and developing robust design conditions using ellipsoidal sets~\cite{peaucelle2005ellipsoidal}.}

%The choice of a diagonal gain is common in the extremum-seeking control system literature.  As a result, it is noted that the system does not converge using the diagonal structure gain. However, with the designed gain in this article, it was possible to ensure the extremum-seeking control system stability with saturated actuators. Also note that the simulation was carried out considering a randomly generated value of $H$ inside the polytope.

%\begin{ack}                               % Place acknowledgements
%Partially supported by the Roman Senate.  % here.
%\end{ack}

%\newpage
\section*{Appendix}

\appendix

\section{\textcolor{black}{Averaging Theory for Lipschitz Continuous Right-Hand Sides}}\label{sec:appendix}

Consider a system of the form
\begin{equation*}
\dot x = \varepsilon f(t,x,\varepsilon), \qquad x(0)=x_0,
\end{equation*}
where $\varepsilon > 0$ is a small parameter, and $f:\mathbb{R}_{+}\times\mathbb{R}^n\times[0,\varepsilon_0]\to\mathbb{R}^n$ is $T$-periodic in time, continuous in $t$, and globally Lipschitz in $x$ with a constant $L>0$, uniformly in $t$. Define the average vector field
\begin{equation*}
\bar f(x) = \frac{1}{T}\int_0^T f(s,x,0)\,ds
\end{equation*}
and consider the corresponding average system
\begin{equation*}
\dot y = \varepsilon \bar f(y), \qquad y(0)=x_0.
\end{equation*}
Even though the function $f$ may not be differentiable with respect to $x$ everywhere, as is the case with functions like the standard saturation, the averaging approach remains valid under the Lipschitz condition. By expressing the solutions of both the original and average systems in integral form,
\begin{equation*}
x(t) \!=\! x_0 \!+\! \varepsilon \!\int_0^t \!\!f(s,x(s),\varepsilon)\,ds, \quad
y(t) \!=\! x_0 \!+\! \varepsilon \!\int_0^t \!\!\bar f(y(s))\,ds,
\end{equation*}
and defining the difference $z(t) := x(t) - y(t)$, one obtains
\begin{equation*}
z(t) = \varepsilon \int_0^t \big[f(s,x(s),\varepsilon) - \bar f(y(s))\big]\, ds.
\end{equation*}
This difference can be split into two terms:
\begin{eqnarray*}
&&f(s,x(s)) - \bar f(y(s)) \nonumber \\
&&= \big[f(s,x(s),\varepsilon) \!-\! f(s,y(s),\varepsilon)\big] \!+\! \big[f(s,y(s),\varepsilon) \!-\! \bar f(y(s))\big].
\end{eqnarray*}
The first term is controlled using the Lipschitz property of $f$, while the second term, corresponding to the oscillatory part, has zero mean over one period and admits a uniformly bounded primitive in time. Let $M>0$ denote a uniform bound on this primitive. Combining both estimates leads to
\begin{equation*}
\|z(t)\| \le \varepsilon L \int_0^t \|z(s)\|\, ds + 2 M \varepsilon,
\end{equation*}
and applying Gr\"onwall's inequality gives
\begin{equation*}
\|x(t) - y(t)\| \le 2 M \varepsilon e^{\varepsilon L t}.
\end{equation*}
Hence, for times $t$ up to $\mathcal{O}(1/\varepsilon)$, the solutions of the original system remain close to those of the average system, with
\begin{equation*}
\|x(t) - y(t)\| = \mathcal{O}(\varepsilon).
\end{equation*}
This result depends only on the Lipschitz continuity of $f$ in $x$, and not on differentiability, making it directly applicable to systems with saturation-type nonlinearities or other piecewise-smooth right-hand sides.

In addition, if the average system $\dot y = \varepsilon \bar f(y)$ has an asymptotically stable equilibrium, then there exists $0<\varepsilon^*<\varepsilon_0$ such that for all $0<\varepsilon<\varepsilon^*$ (sufficiently small) the following inequality is satisfied:
\begin{equation*}
\sup_{t\ge 0} \|x(t) - y(t)\| \leq C \varepsilon, \quad C>0,
\end{equation*}
meaning the approximation remains $\mathcal{O}(\varepsilon)$--close for all times, even when the right-hand side is non-differentiable but Lipschitz (e.g., {saturation or dead-zone functions}). 
%
%the approximation remains valid for all time, even when the right-hand side %is non-differentiable but Lipschitz (e.g., a saturation function).

The results above can be directly derived as a particular case from the more general averaging theorem for systems with discontinuous right-hand sides in \cite{plotnikov1979averaging}.

\bibliographystyle{abbrv}

\bibliography{autosam}      
\end{document}